\documentclass[11pt]{amsart}
\usepackage{amssymb}
\usepackage{amsthm}

\usepackage{amsfonts}
\usepackage{amsmath}
\usepackage[abbrev]{amsrefs}
\usepackage{enumerate}

\input cyracc.def

\newtheorem{Theorem}{Theorem}[section]

\newtheorem{Lemma}[Theorem]{Lemma}
\newtheorem{Corollary}[Theorem]{Corollary}
\newtheorem{Proposition}[Theorem]{Proposition}
\newtheorem{Definition}[Theorem]{Definition}

\theoremstyle{definition}
\newtheorem{Remark}[Theorem]{Remark}

\newcommand{\R}{\mathbb{R}}

\newcommand{\N}{\mathbb{N}}

\newcommand{\A}{\mathcal{H}}
\newcommand{\B}{\hat{\mathcal{H}}}

\DeclareMathOperator{\esssup}{ess\,sup}

\DeclareMathOperator*{\sign}{sign}

\begin{document}

\title{On the dual of Ces{\`a}ro function space.}

\keywords{Ces{\`a}ro function space, dual space, Radon-Nikodym property, slice}
\subjclass[2010]{46E30, 46B20, 46B42, 46B22}

\author{Anna Kami\'{n}ska}
\address{Department of Mathematical Sciences
The University of Memphis, TN 38152-3240}
\email{kaminska@memphis.edu}

\author{Damian Kubiak}
\address{Department of Mathematical Sciences
The University of Memphis, TN 38152-3240}
\email{dmkubiak@memphis.edu}

\date{December 8, 2011}

\begin{abstract}
The goal of this paper is to present an isometric representation of the dual space to Ces{\`a}ro function space $C_{p,w}$, $1<p<\infty$, induced by arbitrary positive weight function $w$ on interval $(0,l)$ where $0<l\leqslant\infty$. For this purpose given a strictly decreasing nonnegative function $\Psi$ on $(0,l)$, the notion of essential $\Psi$-concave majorant $\hat f$ of a measurable function $f$ is introduced and investigated.  As applications it is shown that every slice of the unit ball of the Ces{\`a}ro function space  has diameter 2. Consequently  Ces{\`a}ro function spaces do not have the Radon-Nikodym property, are not locally uniformly convex and they are not dual spaces.
\end{abstract}

\maketitle

\section{Introduction}

Ces{\`a}ro function spaces were studied for the first time in 1970 \cite{MR0276751}. These spaces have been defined analogously to the Ces{\`a}ro sequence spaces that appeared two years earlier in \cite{pr68} when the  Dutch Mathematical Society posted a problem  to find a representation of their dual spaces. This problem was resolved by Jagers \cite{MR0348444} in 1974  who gave an explicit isometric description of the dual of Ces{\`a}ro sequence space. In 1996 different, isomorphic description due to Bennett appeared in \cite{MR1317938}. For a long time Ces{\`a}ro function spaces have not attracted a lot of attention contrary to their sequence counterparts. In fact there is quite rich literature concerning different topics studied in Ces{\`a}ro sequence spaces as for instance in \cites{MR1608225, MR1810056, MR1744077, MR1972393, MR1904283, MR1772119}. However recently in a series of papers \cites{MR2431042,MR2650988,MR2639977}, Astashkin and Maligranda started to study   thoroughly the structure of Ces{\`a}ro function spaces.
Among others, in \cite{MR2639977} they investigated dual spaces for classical Ces{\`a}ro function spaces $Ces_p$ induced by the weight $w(x) = x^{-1}$ for $1<p<\infty$. Their description can be viewed as being analogous to one given for sequence spaces in \cite{MR1317938}. They found a Banach space equipped with a norm equivalent to the dual norm, which is an isomorphic representation of the dual space.

Here, in this paper, we compute precisely the dual norm of the Ces{\`a}ro function space $C_{p,w}$ on $(0,l)$, $0 <l\leqslant \infty$, generated by $1<p<\infty$ and an arbitrary positive weight $w$.
A description presented in this paper resembles the approach of Jagers \cite{MR0348444} for sequence spaces; however, the techniques are more involved due to necessity of dealing with measurable functions instead of sequences.

An introductory section 2 of this paper is devoted to $\Psi$-concave functions and essential $\Psi$-concave majorants of measurable functions and it can of independent interest. In this section $\Psi$ is a nonnegative strictly decreasing  function on the interval $I=(a,b) \subset \R$. The notion related to $\Psi$-concavity was defined by Beckenbach in 1937 \cite[cf. Section 84, p. 240]{MR0442824}. We introduce here also a new notion of essential $\Psi$-concave majorant $\hat f$ of a measurable function $f$, which is a key to study a representation of dual spaces.  We discuss several properties of $\Psi$-concave functions as well as the existence, continuity or differentiability of $\Psi$-concave majorant $\hat f$ of an arbitrary measurable function $f$.

   In main section 3, we give an isometric description of the dual of Ces{\`a}ro function spaces with arbitrary weight function on finite or infinite interval $(0,l)$. We treat finite and infinite case in an unified way, opposite to the isomorphic description given in \cite{MR2639977}. It is also worth to mention that in the process of showing our results, we do not use Hardy inequality at all, an essential tool in studying the space $Ces_p$.

In section 4, applying techniques developed in studying duality we prove that every slice of the unit ball of $C_{p,w}$, $1<p<\infty$, has diameter $2$. Consequently, in a final part we state several corollaries as that Ces{\`a}ro function spaces do not have the Radon-Nikodym property, neither strongly exposed nor denting points, as well as they are not locally uniformly rotund or that they are not dual spaces.

Throughout this paper, the terms decreasing or increasing mean non-increasing or non-decreasing, respectively. By $m$ we denote the Lebesgue measure on the real line  $\R$.
For an interval $I\subset\R$ by $L_{0}(I)$ we denote the set of all (equivalence classes of extended) real valued Lebesgue measurable functions on $I$. The positive cone of $L_0(I)$ is denoted $L_0^{+}(I)=\{f\in L_0(I):f\geqslant 0\ a.e.\}$.

For a Banach space $(X,\|\cdot\|)$ by $B_X$ and $S_X$ we denote the unit ball and the unit sphere of $X$, and by $X^*$ the dual space of $X$. Any Banach space $E=E(I)\subset L_0(I)$ with norm $\|\cdot\|$ satisfying the condition that $f\in E$ and $\|f\|\leqslant\|g\|$  whenever $0\leqslant f\leqslant g$ a.e., $f\in L_0(I)$ and $g\in E$, is called a \emph{Banach function lattice}.
An element $f$ in a Banach function lattice $E$ is called \emph{order continuous} if for every $0\leqslant f_n\leqslant|f|$ a.e. such that $f_n\downarrow 0$ a.e. it holds $\|f_n\|\downarrow 0$. We say that $E$ is \emph{order continuous} if every element in $E$ is order continuous. A Banach function lattice $(E,\|\cdot\|)$ has the \emph{Fatou property} if for any sequence $(f_n)\subset E$ and any $f\in L_0(I)$ such that $0\leqslant f_n\leqslant f$ a.e., $f_n\uparrow f$ a.e. and $\sup_n\|f_n\|<\infty$ it holds $f\in E$ and $\|f\|=\lim_n\|f_n\|$.

\section{$\Psi$-concave functions and essential $\Psi$-concave majorants}

In this section we fix $I=(a,b)\subset\R$ to be an open (finite or infinite) interval and $\Psi:I\to \R_{+}$ to be a strictly decreasing function on $I$.
We first collect a number of properties of $\Psi$-concave functions. Some of them are certainly known \cite[cf Section 84, p. 240]{MR0442824} but we provide their proofs here for the sake of completeness. Next we introduce a notion of essential $\Psi$-concave majorant $\hat f$ of $f\in L_0(I)$ and discuss a number of its properties like existence, continuity and differentiability.
This section can be of independent interest.

Recall a definition of $\Psi$-concave function \cite{MR0348444}.
\begin{Definition}
A function $f:I\to \R$ is called \emph {$\Psi$-concave} (respectively \emph{strictly $\Psi$-concave}) on $I$ if for all $x<y<z$ in $I$,
\begin{equation}
\label{psidef}
\left|
\begin{array}{ccc}
   1 & 1 & 1\\
   \Psi(x)& \Psi(y)& \Psi(z) \\
    f(x)& f(y) & f(z)
\end{array} \right| \geqslant 0\quad\text{(respectively, }>0\text{).}
\end{equation}
\end{Definition}

It is easy to check that inequality (\ref{psidef}) is equivalent to
\begin{equation}
\label{psiconcave0}
\frac{f(x)-f(y)}{\Psi(x)-\Psi(y)}\leqslant\frac{f(y)-f(z)}{\Psi(y)-\Psi(z)}\quad\text{(respectively }<\text{) for all }x<y<z\text{ in }I\text{.}
\end{equation}
It is also possible to rewrite (\ref{psiconcave0}) as
\begin{equation}
\label{psiconcave3}
\frac{f(x)-f(y)}{\Psi(x)-\Psi(y)}\leqslant\frac{f(x)-f(z)}{\Psi(x)-\Psi(z)}\leqslant\frac{f(y)-f(z)}{\Psi(y)-\Psi(z)}\quad\text{(respectively }<\text{) for all } x<y<z\text{ in }I\text{.}
\end{equation}

If the interval $I=(a,b)$ is finite and $\Psi(x)=b-x$, $x\in I$, then $\Psi$-concavity on $I$ is just usual concavity.

The following definition will be also useful.
\begin{Definition}
We say that a function $f:I\mapsto \R$ is \emph{$\Psi$-affine} on $I$, if $f(x)=A\Psi(x)+B$, $x\in I$, for some constants $A$ and $B$.
\end{Definition}

For arbitrary interval $J$ we say that $f:J\to\R$ is \emph{$\Psi$-concave on $J$} if it is $\Psi$-concave on the interior of $J$. Similarly in the case of $\Psi$-affine function.

Now, similarly as done for example in \cite{MR0442824} in context of convex functions, we show basic properties of $\Psi$-concave functions. Let $f:I\to\R$ be $\Psi$-concave on $I$. Define for $x\in I$,
\[
D_{\Psi}^{+}f(x)=\lim_{y\to x^{+}}\frac{f(y)-f(x)}{\Psi(y)-\Psi(x)}\ \ \text{ and }\ \ D_{\Psi}^{-}f(x)=\lim_{y\to x^{-}}\frac{f(y)-f(x)}{\Psi(y)-\Psi(x)}\text{.}
\]
In order to see the existence and finiteness of the above quantities, it is enough to observe that for any $w<x<y<z<u$ in $I$, by (\ref{psiconcave3}) it follows that
\[
\frac{f(w)-f(y)}{\Psi(w)-\Psi(y)}\leqslant\frac{f(x)-f(y)}{\Psi(x)-\Psi(y)}\leqslant\frac{f(y)-f(z)}{\Psi(y)-\Psi(z)}\leqslant\frac{f(y)-f(u)}{\Psi(y)-\Psi(u)}\text{,}
\]
and hence the left side of the inequality
\[
\frac{f(x)-f(y)}{\Psi(x)-\Psi(y)}\leqslant\frac{f(y)-f(z)}{\Psi(y)-\Psi(z)}
\]
increases as $x\uparrow y$ and the right side decreases as $z\downarrow y$. It follows that $D_{\Psi}^{-}f(y)$, $D_{\Psi}^{+}f(y)$ exist and $D_{\Psi}^{-}f(y)\leqslant D_{\Psi}^{+}f(y)$ for all $y\in I$. Monotonicity of $D_{\Psi}^{+}f$ and $D_{\Psi}^{-}f$ follows again from (\ref{psiconcave3}), namely, for any $w<x<y<z$ in $I$,
\[
D_{\Psi}^{+}f(w)=\lim_{y\downarrow w}\frac{f(y)-f(w)}{\Psi(y)-\Psi(w)}\leqslant\frac{f(w)-f(x)}{\Psi(w)-\Psi(x)}\leqslant\frac{f(y)-f(z)}{\Psi(y)-\Psi(z)}\leqslant\lim_{y\uparrow z}\frac{f(y)-f(z)}{\Psi(y)-\Psi(z)}=D_{\Psi}^{-}f(z)\text{.}
\]
Hence $D_{\Psi}^{-}f(w)\leqslant D_{\Psi}^{+}f(w)\leqslant D_{\Psi}^{-}f(z)\leqslant D_{\Psi}^{+}f(z)$ for all $w,z\in I$ such that $w<z$.

In fact $D_{\Psi}^{+}f$ is right continuous if $\Psi$ is right continuous. Indeed, by monotonicity of $D_{\Psi}^{+}f$ we have that $\lim_{x\to w^{+}}D_{\Psi}^{+}f(x)$ exists for any $w\in I$. Since for any $y>x$,
\[
D_{\Psi}^{+}f(x)\leqslant\frac{f(y)-f(x)}{\Psi(y)-\Psi(x)}\text{,}
\]
and since $f$ and $\Psi$ are right continuous, $\lim_{x\to w^{+}}D_{\Psi}^{+}f(x)\leqslant\frac{f(y)-f(w)}{\Psi(y)-\Psi(w)}$, $y>x>w$. It follows that
\[
\lim_{x\downarrow w}D_{\Psi}^{+}f(x)\leqslant\lim_{y\downarrow w}\frac{f(y)-f(w)}{\Psi(y)-\Psi(w)}=D_{\Psi}^{+}f(w)\text{.}
\]
On the other hand, we know that $D_{\Psi}^{+}f(w)\leqslant D_{\Psi}^{+}f(x)$ for all $w<x$, and so $\lim_{x\downarrow w}D_{\Psi}^{+}f(x)=D_{\Psi}^{+}f(w)$ for all $w\in I$. Similarly one can show the left continuity of $D_{\Psi}^{-}f$ under assumption of left continuity on $\Psi$.

The next proposition summarizes our discussion as follows.
\begin{Proposition}
\label{Dpsi}
If $f$ is $\Psi$-concave on $I$ then $D_{\Psi}^{+}f(x)$, $D_{\Psi}^{-}f(x)$ exist, are finite and $D_{\Psi}^{-}f(x)\leqslant D_{\Psi}^{+}f(x)$ for all $x\in I$. Moreover,  $D_{\Psi}^{+}f$, $D_{\Psi}^{-}f$ are increasing functions on $I$. If $\Psi$ is right-continuous on $I$ then so is $D_{\Psi}^{+}f$. Similarly, if $\Psi$ is left-continuous on $I$ then so is $D_{\Psi}^{-}f$. For any fixed $y\in I$ the ratio $\frac{f(x)-f(y)}{\Psi(x)-\Psi(y)}$ increases as $x\uparrow y$ and the ratio $\frac{f(y)-f(z)}{\Psi(y)-\Psi(z)}$ decreases as $z\downarrow y$.
\end{Proposition}

The following basic fact will be also useful.
\begin{Lemma}
\label{continuity}
If $\Psi$ is right-, left-, absolutely or Lipschitz continuous then any $\Psi$-concave function $f$  on $I$ has the same property.
\begin{proof}
Let $[c,d]\subset I$ and $a<c_1<c$ and $d<d_1<b$, by $\Psi$-concavity of $f$ we get
\[
\frac{f(c_1)-f(c)}{\Psi(c_1)-\Psi(c)}\leqslant\frac{f(x)-f(y)}{\Psi(x)-\Psi(y)}\leqslant\frac{f(d)-f(d_1)}{\Psi(d)-\Psi(d_1)}\quad\text{for all }c\leqslant x<y\leqslant d\text{,}
\]
hence
\[
\left|\frac{f(x)-f(y)}{\Psi(x)-\Psi(y)}\right|\leqslant\max\left\{\left|\frac{f(c_1)-f(c)}{\Psi(c_1)-\Psi(c)}\right|,\left|\frac{f(d)-f(d_1)}{\Psi(d)-\Psi(d_1)}\right|\right\}\text{.}
\]
Denoting the right hand side by $K$ we get that
\begin{equation*}
\label{Lipschitz}
|f(x)-f(y)|\leqslant K|\Psi(x)-\Psi(y)|\quad\text{for all }x,y\in[c,d]\text{.}
\end{equation*}
The claim follows. 
\end{proof}
\end{Lemma}

\begin{Lemma}
\label{monotone}
Let $f$ be $\Psi$-concave on $I$. If $D_{\Psi}^{+}f(x)\geqslant 0$, $x\in I$, then $f$ is decreasing on $I$. If $D_{\Psi}^{+}f(x)\leqslant 0$, $x\in I$, then $f$ is increasing on $I$.
\begin{proof}
Let $D_{\Psi}^{+}f(x)\geqslant 0$. Since the ratio $\frac{f(z)-f(x)}{\Psi(z)-\Psi(x)}$ decreases as $z\downarrow x$, it follows that $\frac{f(z)-f(x)}{\Psi(z)-\Psi(x)}\geqslant 0$ for $z>x$, and hence $f(z)\leqslant f(x)$ by monotonicity of $\Psi$. Since $x\in I$ is arbitrary, we get that $f$ is decreasing. The proof of another case is similar.
\end{proof}
\end{Lemma}

\begin{Lemma}
\label{Psiinfty}
Let function $f\geqslant 0$ be $\Psi$-concave on $I$. If $\lim_{x\to a^{+}}\Psi(x)=\infty$ then $\lim_{x\to a^{+}}D_{\Psi}^{+}f(x)\geqslant 0$.
\begin{proof}
Since function $D_{\Psi}^{+}f$ is increasing, $\lim_{x\to a^{+}}D_{\Psi}^{+}f(x)$ exists or is equal to $-\infty$. Suppose that $\lim_{x\to a^{+}}D_{\Psi}^{+}f(x)=C<0$. It follows that there exists $x_0>a$ such that $-\infty<D:=D_{\Psi}^{+}f(x_0)<0$, so $D_{\Psi}^{-}f(x)\leqslant D$ for $x\in(a,x_0)$. It follows that for all $z<x$, $\frac{f(x)-f(z)}{\Psi(x)-\Psi(z)}\leqslant D<0$, which gives $f(z)\leqslant D\Psi(z)-D\Psi(x)+f(x)$. Now, keeping $x\in(a,x_0)$ fixed and taking $z\to a^{+}$ we would get that $f(z)<0$ for $z$ close enough to $a$, which contradicts the condition $f\geqslant 0$.
\end{proof}
\end{Lemma}
Lemmas \ref{monotone} and \ref{Psiinfty} imply the following corollary.
\begin{Corollary}
\label{fdecr}
If a function $f\geqslant 0$ is $\Psi$-concave on $I$ and $\lim_{x\to a^{+}}\Psi(x)=\infty$ then $f$ is decreasing on $I$.
\end{Corollary}

Observe that inequality (\ref{psiconcave0}) can also be equivalently written as
\begin{equation}
\label{psiconcave1}
f(y)\geqslant\frac{\Psi(y)-\Psi(z)}{\Psi(x)-\Psi(z)}f(x)+\frac{\Psi(x)-\Psi(y)}{\Psi(x)-\Psi(z)}f(z)\text{ for all }x,y,z\in I\text{ with }x<y<z\text{.}
\end{equation}

\begin{Lemma}
\label{support}
A function $f:I\to\R$ is $\Psi$-concave on $I$ if and only if for each $y\in I$ there is at least one function $T(x)=f(y)+A(\Psi(x)-\Psi(y))$ such that $A\in[D_{\Psi}^{-}f(y),D_{\Psi}^{+}f(y)]$ and $f(x)\leqslant T(x)$ for $x\in I$.
\begin{proof}
If $f$ is $\Psi$-concave on $I$ and $y\in I$ then for any $A\in[D_{\Psi}^{-}f(y),D_{\Psi}^{+}f(y)]$, 
\[
\frac{f(x)-f(y)}{\Psi(x)-\Psi(y)}\geqslant A\quad\text{or }\leqslant A\text{,}
\]
if $x>y$ or $x<y$, respectively. In any case $f(x)\leqslant A\Psi(x)+f(y)-A\Psi(y)=T(x)$ for all $x\in I$.

Conversely, suppose that for each $y\in I$ there is at least one function $T(x)=f(y)+A(\Psi(x)-\Psi(y))$ such that $f(x)\leqslant T(x)$ for $x\in I$. Let $x, y, z\in I$ be such that $x<y<z$. Denoting $\alpha=\frac{\Psi(y)-\Psi(z)}{\Psi(x)-\Psi(z)}$ we get $\Psi(y)=\alpha\Psi(x)+(1-\alpha)\Psi(z)$, $\alpha\in[0,1]$. It follows that $f(y)=T(y)=\alpha T(x)+(1-\alpha)T(z)\geqslant\alpha f(x)+(1-\alpha)f(z)$. Hence,  in a view of (\ref{psiconcave1}), $f$ is $\Psi$-concave.
\end{proof}
\end{Lemma}

The following lemma will be useful.
\begin{Lemma}
\label{psiaff1}
\label{psiaff2}
\label{strpsiconc}

Let $f$ be $\Psi$-concave on $I$. The function $f$ is strictly $\Psi$-concave on $I$ if and only if there is no interval $(c,d)\subset I$ on which $f$ is $\Psi$-affine. The function $f$ is $\Psi$-affine on $I$ if and only if $D_{\Psi}^{+}f$ is constant on $I$. The function $f$ is strictly $\Psi$-concave on $I$ if and only if $D_{\Psi}^{+}f$ is strictly increasing on $I$.
\begin{proof}
We prove only the first part. The proof of the second part is similar.

If there exists an interval $(c,d)\subset I$ on which $f(x)=A\Psi(x)+B$, then $\frac{f(x)-f(y)}{\Psi(x)-\Psi(y)}=A$ for all $x,y\in(c,d)$. Hence $f$ is not strictly $\Psi$-concave on $I$.
Conversely, if $f$ is not strictly $\Psi$-concave on $I$, then from (\ref{psiconcave0}) we get that there exist $w<x<u$ in $I$ such that
\[
A:=\frac{f(w)-f(x)}{\Psi(w)-\Psi(x)}=\frac{f(x)-f(u)}{\Psi(x)-\Psi(u)}\text{.}
\]
It follows from (\ref{psiconcave3}) that for any $y$ in $(x,u)$,
\[
A=\frac{f(w)-f(x)}{\Psi(w)-\Psi(x)}\leqslant\frac{f(w)-f(y)}{\Psi(w)-\Psi(y)}\leqslant\frac{f(w)-f(u)}{\Psi(w)-\Psi(u)}\leqslant\frac{f(x)-f(u)}{\Psi(x)-\Psi(u)}=A\text{.}
\]
Hence $f(y)=A\Psi(y)+f(w)-A\Psi(w)$ for all $y\in(x,u)$, that is $f$ is $\Psi$-affine on $(x,u)$, a contradiction.
\end{proof}
\end{Lemma}


 In the case when $\Psi$ is a continuous function, we can define $\Psi$-concavity in one more equivalent way. Namely, denoting $\alpha=\frac{\Psi(y)-\Psi(z)}{\Psi(x)-\Psi(z)}$ for $x<y<z$ in $I$, we get that $\Psi(y)=\alpha\Psi(x)+(1-\alpha)\Psi(z)$, and so (\ref{psiconcave1}) can be written as
\begin{equation}
\label{psiconcave2}
f(\Psi^{-1}(\alpha\Psi(x)+(1-\alpha)\Psi(z)))\geqslant\alpha f(x)+(1-\alpha)f(z)\text{.}
\end{equation}
If $\Psi$ is continuous then for any $x,z\in I$, say $x\leqslant z$, and any $\alpha\in[0,1]$ there exists $y\in[x,z]$ such that $\Psi(y)=\alpha\Psi(x)+(1-\alpha)\Psi(z)$. Then inequality (\ref{psiconcave2}) holds true for all $\alpha\in[0,1]$ and all $x,z\in I$, that is function $f\circ\Psi^{-1}$ is concave on $\Psi(I)$. Furthermore, in this case, by induction it can be shown that
\begin{equation}
\label{psiconcave}
f\left(\Psi^{-1}\left(\sum_{i=1}^{n}\alpha_i\Psi(y_i)\right)\right)\geqslant\sum_{i=1}^{n}\alpha_i f(y_i)\quad\text{for all }\alpha_i\geqslant 0\text{, }\sum_{i=1}^{n}\alpha_i=1\text{ and }(y_1, y_2,\ldots, y_n)\in I^n\text{.}
\end{equation}

We have the following lemma.
\begin{Lemma}
\label{concaveEqv}
If $\Psi$ is a continuous function on $I$ then $f$ is $\Psi$-concave on $I$ if and only if $f\circ\Psi^{-1}$ is concave on $\Psi(I)$.
\begin{proof}
Only one direction requires proof. Suppose that $f\circ\Psi^{-1}$ is concave on $\Psi(I)$. Let $x,y,z\in I$ be arbitrary, $x<y<z$ and $u=\Psi(x)$, $v=\Psi(z)$. Since $\Psi$ is one-to-one, $x=\Psi^{-1}(u)$ and $z=\Psi^{-1}(v)$. By concavity of $f\circ\Psi^{-1}$ we get that $(f\circ\Psi^{-1})(\alpha u+(1-\alpha) v)\geqslant \alpha(f\circ\Psi^{-1})(u)+(1-\alpha)(f\circ\Psi^{-1})(v)$ for all $\alpha\in[0,1]$. It follows that $f(\Psi^{-1}(\alpha\Psi(x)+(1-\alpha)\Psi(z)))\geqslant \alpha f(x)+(1-\alpha)f(z)$. Since $y\in(x,z)$ there exists $\alpha\in[0,1]$ such that $\Psi(y)=\alpha\Psi(x)+(1-\alpha)\Psi(z)$. This gives inequality  (\ref{psiconcave1}) and so $f$ is $\Psi$-concave on $I$.
\end{proof}
\end{Lemma}

The following notion of essential $\Psi$-concave majorant is crucial for characterization of the dual space to Ces{\`a}ro function spaces.
\begin{Definition}
For any function $f\in L_0^{+}(I)$ we define its \emph{essential $\Psi$-concave majorant} $\hat{f}$ by
\begin{equation*}
\begin{split}
\hat{f}(y)
:=\inf\bigg\{M>0:m^{(n)}\Big\{(y_1,\ldots,y_n)\in I^n:&\sum_{i=1}^{n}\alpha_if(y_i)>M, \sum_{i=1}^{n}\alpha_i=1, \alpha_i\geqslant0, i=1,\ldots, n, \\
& \Psi(y)=\sum_{i=1}^{n}\alpha_i\Psi(y_i)\Big\}=0\text{, }n\in\N\bigg\}\text{,}\quad y\in I\text{,}
\end{split}
\end{equation*}
where $m^{(n)}$ is the Lebesgue product measure on $I^n$. For arbitrary function $f\in L_0(I)$ we define $\hat{f}=\widehat{|f|}$.
\end{Definition}

The above definition should be compared to one of concave majorants given in 1970 by Peetre \cite{MR0272960}.

The remaining results of this section describe several properties of
$\hat{f}$.
First we give conditions on $f$ under which the essential $\Psi$-concave majorant $\hat{f}$ is finite on $I$.

\begin{Lemma}
\label{finitehat1}
Let $f\in L_0^{+}(I)$. If $\esssup_{x\in(y,b)}f(x)<\infty$ and $\esssup_{x\in(a,y)}\frac{f(x)}{\Psi(x)}<\infty$ for all $y\in I$ then $\hat{f}<\infty$ on $I$.
\begin{proof}
Let $f\in L_0^{+}(I)$ and $y\in I$. Suppose that
$$
A_y:=\esssup_{x\in(y,b)}f(x)<\infty\ \ \ \text{and}\ \ \ B_y:=\esssup_{x\in(a,y)}\frac{f(x)}{\Psi(x)}<\infty.
$$
For any $n\in\N$, if $\Psi(y)=\sum_{i=1}^{n}\alpha_i\Psi(y_i)$, $\sum_{i=1}^{n}\alpha_i=1$ we have that
\begin{equation*}
\begin{split}
\sum_{i=1}^{n}\alpha_if(y_i)&=\sum_{y_i<y}\alpha_if(y_i)+\sum_{y_i\geqslant y}\alpha_if(y_i)=\sum_{y_i<y}\alpha_i\Psi(y_i)f(y_i)/\Psi(y_i)+\sum_{y_i\geqslant y}\alpha_if(y_i)\leqslant B_y\Psi(y)+A_y
\end{split}
\end{equation*}
except possibly some subset of the set
\[
C:=\bigcup\Big\{(y_1,\ldots,y_n)\in I^n:\max_{i\in I_1:y_i<y}(f(y_i)/\Psi(y_i))>B_y\Big\}\cup\Big\{(y_1,\ldots,y_n)\in I^n:\max_{i\in I_2:y_i\geqslant y}(f(y_i))>A_y\Big\}\text{,}
\]
where the union is taken over all partitions of $\{1,2,\ldots,n\}$ into two disjoint nonempty sets $I_1$, $I_2$.
It is not difficult to see that $m^{(n)}C=0$, whence it follows that $\hat{f}(y)\leqslant B_y\Psi(y)+A_y<\infty$.
\end{proof}
\end{Lemma}

\begin{Lemma}
\label{psiconcavehat}
Let $f\in L_0^{+}(I)$. If $\hat{f}<\infty$ on $I$ then $\hat{f}$ is $\Psi$-concave on $I$.
\begin{proof}
The proof is similar to one for concave majorants \cite[p. 47]{KPS}.

Let $x<y<z$ in $I$. We will show inequality (\ref{psiconcave1}) for $\hat{f}$. Let $\alpha=\frac{\Psi(y)-\Psi(z)}{\Psi(x)-\Psi(z)}$ and $\epsilon>0$ be arbitrary. From the definition of $\hat{f}(x)$ it follows that there exist $j\in\N$ and a set
\[
 B=\left\{(x_1^{\epsilon},x_2^{\epsilon}, \ldots, x_j^{\epsilon} )\in I^j:\sum_{i=1}^{j}\alpha_i^{\epsilon}f(x_i^{\epsilon})>\hat{f}(x)-\epsilon/2, \Psi(x)=\sum_{i=1}^{j}\alpha_i^{\epsilon}\Psi(x_i^{\epsilon}), \sum_{i=1}^{j}\alpha_i^{\epsilon}=1\right\}
\]
with $m^{(j)}B>0$. Denoting $\alpha_i^{'\epsilon}=\alpha\alpha_i^{\epsilon}$, $i=1, 2, \ldots, j$, we get that for all $(x_1^{\epsilon},x_2^{\epsilon}, \ldots, x_j^{\epsilon})\in B$, $\alpha\hat{f}(x)\leqslant\sum_{i=1}^j\alpha_i^{'\epsilon}f(x_i^{\epsilon})+\epsilon/2$, where $\sum_{i=1}^{j}\alpha_i^{'\epsilon}=\alpha$, $\alpha_i^{'\epsilon}>0$, $\alpha\Psi(x)=\sum_{i=1}^j\alpha_i^{'\epsilon}\Psi(x_i^{\epsilon})$. Similarly, by definition of $\hat{f}(z)$, there exist $k\in\N$ and a set
\[
 C=\left\{(z_1^{\epsilon},z_2^{\epsilon}, \ldots, z_k^{\epsilon} )\in I^k:\sum_{i=1}^{k}\beta_i^{\epsilon}f(z_i^{\epsilon})>\hat{f}(z)-\epsilon/2, \Psi(z)=\sum_{i=1}^{k}\beta_i^{\epsilon}\Psi(z_i^{\epsilon}), \sum_{i=1}^{k}\beta_i^{\epsilon}=1\right\}
\]
 with $m^{(k)}C>0$. Denoting $\beta_i^{'\epsilon}=(1-\alpha)\beta_i^{\epsilon}$, $i=1, 2, \ldots, k$, we get that for all $(z_1^{\epsilon},z_2^{\epsilon}, \ldots, z_k^{\epsilon})\in C$, $(1-\alpha)\hat{f}(z)\leqslant\sum_{i=1}^k\beta_i^{'\epsilon}f(z_i^{\epsilon})+\epsilon/2$, where $\sum_{i=1}^{k}\beta_i^{'\epsilon}=1-\alpha$, $\beta_i^{'\epsilon}>0$, $(1-\alpha)\Psi(z)=\sum_{i=1}^k\beta_i^{'\epsilon}\Psi(z_i^{\epsilon})$. Denoting now $\gamma_i^{\epsilon}=\alpha_i^{'\epsilon}$, $y_i^{\epsilon}=x_i^{\epsilon}$ for $i=1, 2, \ldots, j$, $\gamma_{i+j}^{\epsilon}=\beta_i^{'\epsilon}$, $y_{i+j}^{\epsilon}=z_i^{\epsilon}$ for $i = 1, 2, \ldots k$, and $n=j+k$ we get that $\sum_{i=1}^n\gamma_i^{\epsilon}=1$ and $\sum_{i=1}^n\gamma_i^{\epsilon}\Psi(y_i^{\epsilon})=\alpha\Psi(x)+(1-\alpha)\Psi(z)=\Psi(y)$. It follows that $\alpha\hat{f}(x)+(1-\alpha)\hat{f}(z)\leqslant\sum_{i=1}^n\gamma_i^{\epsilon}f(y_i^{\epsilon})+\epsilon\leqslant\hat{f}(y)+\epsilon$ a.e. Since $\epsilon$ was arbitrary the claim follows.
\end{proof}
\end{Lemma}

Recall that if $C$ is a measurable subset  of $\R$ and $y\in\R$ then $y$ is called a \emph{point of density} of $C$ if
\[
\lim_{\substack{m(x,z)\to 0\\ y\in (x,z)}}\frac{m(C\cap(x,z))}{m(x,z)}=1\text{.}
\]
It is known that if $C$ is a measurable subset of $\R$ then almost every $x\in C$ is a point of density of $C$ \cite[p. 106]{MR2129625} \cite[p. 141]{MR924157}.

\begin{Lemma}
\label{hatgr}
If $\Psi$ is a continuous function and $\lim_{x\to a^{+}}\Psi(x)=\infty$ then for any function $f\in L_0^{+}(I)$ with $\hat{f}<\infty$ on $I$, it holds that $f\leqslant\hat{f}$ a.e. on $I$, and $\hat{f}$ is also continuous on $I$.
\begin{proof}
Suppose there exist $\epsilon>0$ and a set $C\subset I$ with $mC>0$ such that $f\geqslant\hat{f}+\epsilon$ on $C$. Without loss of generality we assume that all points in $C$ are points of density of $C$. It follows that for all $y\in C$ and all $x<y<z$ in $I$, $m(C\cap(x,z))>0$.

First we show that for any $x<z$ such that $m(C\cap(x,z))>0$ there is $y\in C\cap(x,z)$ for which
$m(C\cap(x,y))>0$ and $m(C\cap(y,z))>0$.
Let $c=\sup_{y\in C\cap(x,z)}m(C\cap(x,y))=0$ and $d=\inf_{y\in C\cap(x,z)}m(C\cap(y,z))=0$. It follows that $c<d$, 
which gives $m(C\cap(x,y))>0$ and $m(C\cap(y,z))>0$ for all $c<y<d$.

Since $\Psi$ is a continuous function, by Lemmas \ref{psiconcavehat} and \ref{continuity} we get that $\hat{f}$ is continuous. Let $x<z$ in $I$ be such that
\begin{equation}
\label{eps2}
\hat{f}(x)-\hat{f}(z)<\epsilon/2\ \ \text{ and }\ \ m(C\cap(x,z))>0\text{.}
\end{equation}
 By the above there is $y\in C\cap(x,z)$ such that $m(C\cap(x,y))>0$ and $m(C\cap(y,z))>0$.
Consider the set $B:=\{(y_1,y_2)\in (C\cap(x,z))\times (C\cap(x,z)):\Psi(y)=\alpha\Psi(y_1)+(1-\alpha)\Psi(y_2)\text{ for some }\alpha\in(0,1)\}$. It is clear that $B=\{(y_1,y_2)\in (C\cap(x,z))\times (C\cap(x,z)):y_1<y<y_2\quad\text{or }y_2<y<y_1\}=((C\cap(x,y))\times (C\cap(y,z)))\cup((C\cap(y,z))\times (C\cap(x,y)))$ and hence $m^{(2)}B>0$. Observe that by (\ref{eps2}) we have $|\hat{f}(y_1)-\hat{f}(y_2)|<\epsilon/2$ for all $y_1,y_2\in I$ such that $(y_1,y_2)\in B$. Now, for almost all  $(y_1,y_2)\in B$, since $\hat{f}$ is decreasing by Corollary \ref{fdecr}, we have that
\[
\hat{f}(y)\geqslant\alpha f(y_1)+(1-\alpha)f(y_2)\geqslant\alpha \hat{f}(y_1)+(1-\alpha)\hat{f}(y_2)+\epsilon\geqslant\hat{f}(y)+\epsilon/2\text{.}
\]
This is impossible, hence $f\leqslant\hat{f}$ a.e. on $I$.

\end{proof}
\end{Lemma}

\begin{Remark}
\label{remark}
\begin{itemize}
\item[$(1)$] If $f, g\in L_0^{+}(I)$ and $f\leqslant g$ a.e. on $I$ then $\hat{f}\leqslant\hat{g}$. In fact $\sum_{i=1}^n\alpha_if(y_i)\leqslant\sum_{i=1}^n\alpha_ig(y_i)$ for all $(y_1,y_2, \ldots, y_n)\in I^n$ except possibly some set of measure $0$.
\item[$(2)$] Let $\Psi$ be a continuous function on $I$ and $\lim_{x\to a^{+}}\Psi(x)=\infty$. If $f$ is $\Psi$-concave on $I$ then $f=\hat{f}$. Consequently $\hat{\hat{f}}=\hat{f}$ for any function $f\in L_0(I)$ with $\hat{f}<\infty$ on $I$. Indeed, from (\ref{psiconcave}) it follows that $\sum_{i=1}^n\alpha_if(y_i)\leqslant f(y)$ whenever $\Psi(y)=\sum_{i=1}^{n}\alpha_i\Psi(y_i)$ and hence $\hat{f}\leqslant f$ which together with Lemma \ref{hatgr} gives that $f=\hat{f}$.
 \end{itemize}


\end{Remark}

\begin{Lemma}
\label{affine_on_interval}
Let $\Psi$ be a continuous function on $I$ and $\lim_{x\to a^{+}}\Psi(x)=\infty$. Let $f\in L_0^{+}(I)$ be such that $\hat{f}<\infty$ on $I$, $\epsilon>0$ be fixed, $A=\{x\in I:f(x)\geqslant\hat{f}(x)-\epsilon\}$ and $(u,v)\subset I$ be a finite open interval. If $m(A\cap(u,v))=0$ then $\hat{f}$ is $\Psi$-affine on $(u,v)$.
\begin{proof}
Let $y\in (u,v)$ be fixed.
For any $\eta>0$, by definition of $D_{\Psi}^{-}\hat{f}(y)$,
\[
D_{\Psi}^{-}\hat{f}(y)\geqslant \frac{\hat{f}(c)-\hat{f}(y)}{\Psi(c)-\Psi(y)}\geqslant D_{\Psi}^{-}\hat{f}(y)-\eta
\]
for all $c<y$ close enough to $y$. Moreover, the ratio $\frac{\hat{f}(c)-\hat{f}(x)}{\Psi(c)-\Psi(x)}$ is a continuous function of $x$ and by Proposition \ref{Dpsi} it decreases as $x\downarrow y$. Hence for every $\eta>0$ and every $c<y$ there exists $d>y$ arbitrary close to $y$ such that
\[
D_{\Psi}^{-}\hat{f}(y)+\eta\geqslant \frac{\hat{f}(c)-\hat{f}(d)}{\Psi(c)-\Psi(d)}\geqslant D_{\Psi}^{-}\hat{f}(y)-\eta\text{.}
\]

By the above, we construct an increasing sequence $(a_n)\subset(u,y)$ and a sequence  $(b_n)\subset(y,v)$ such that $a_n\to y$, $b_n\to y$ and
\[
s_n:=\frac{\hat{f}(a_n)-\hat{f}(b_n)}{\Psi(a_n)-\Psi(b_n)}\to C:=D_{\Psi}^{-}\hat{f}(y)\quad\text{as }n\to\infty\text{.}
\]
Consider the sequence of functions
\[
S_n(x)=s_n\Psi(x)+\hat{f}(a_n)-s_n\Psi(a_n)\text{, }\ \ x\in I\text{.}
\]
It is clear that $S_n(a_n)=\hat{f}(a_n)$, $S_n(b_n)=\hat{f}(b_n)$ and by (\ref{psiconcave1}), $S_n(x)\leqslant \hat{f}(x)$ for all $x\in(a_n,b_n)$, $n\in\N$.
By Lemma \ref{support}, $\hat{f}(x)\leqslant C\Psi(x)+(\hat{f}(y)-C\Psi(y))$. Hence
\begin{equation*}
\begin{split}
\hat{f}(x)-S_n(x)&=\hat{f}(x)-(s_n\Psi(x)+(\hat{f}(a_n)-s_n\Psi(a_n)))\\
&\leqslant C\Psi(x)+\hat{f}(y)-C\Psi(y)-s_n\Psi(x)-\hat{f}(a_n)+s_n\Psi(a_n)\\
&=(C-s_n)\Psi(x)+\hat{f}(y)-\hat{f}(a_n)+s_n\Psi(a_n)-C\Psi(y)\\
&\leqslant |C-s_n|\Psi(a_1)+\hat{f}(y)-\hat{f}(a_n)+s_n\Psi(a_n)-C\Psi(y)\text{.}
\end{split}
\end{equation*}
It follows that for every $\delta>0$ there exists $N_{\delta}\in\N$ such that for all $n>N_{\delta}$ and for all $x\in (a_n,b_n)$,
\begin{equation*}
0\leqslant \hat{f}(x)-S_n(x)\leqslant\delta\text{.}
\end{equation*}

By the above for all $y\in(u,v)$ there exist $c_y,d_y\in (u,v)$, $c_y<y<d_y$, such that $\hat{f}(x)\geqslant D\Psi(x)+B$ and $\hat{f}(x)-(D\Psi(x)+B)\leqslant\epsilon$ for all $x\in(c_y,d_y)$ where 
$D=(\hat{f}(c_y)-\hat{f}(d_y))/(\Psi(c_y)-\Psi(d_y))$ and $B=(\Psi(c_y)\hat{f}(d_y)-\Psi(d_y)\hat{f}(c_y))/(\Psi(c_y)-\Psi(d_y))$. By definition of the set $A$ it follows that $f(x)\leqslant D\Psi(x)+B$ a.e. on $(c_y,d_y)$. Since function  $g(t)=\hat{f}\chi_{(c_y,d_y)^c}(t)+(D\Psi(t)+B)\chi_{(c_y,d_y)}(t)$, $t\in I$, is $\Psi$-concave on $I$ it must be $\hat{f}=g$ by Remark \ref{remark}. But $g$ is $\Psi$-affine on $(c_y,d_y)$, so
$D_{\Psi}^{+}\hat{f}$ is constant there by Lemma \ref{psiaff2}.

The family of set $(c_y,d_y)$, $y\in(u,v)$, covers each closed subinterval $[u+\epsilon,v-\epsilon]$, $\epsilon>0$. Using compactness we conclude that $D_{\Psi}^{+}\hat{f}$ is constant on $(u,v)$ and by Lemma \ref{psiaff2}, $\hat{f}$ is $\Psi$-affine on $(u,v)$.
\end{proof}
\end{Lemma}

Recall the following theorem concerning convex functions \cite[Corollary 1.3.8, p.23]{MR2178902}.
\begin{Theorem}
\label{convex_conv}
If $f_n:I\to\R$ is a pointwise converging sequence of convex functions, then the limit is also convex. Moreover, the convergence is uniform on any compact subinterval included in the interior of $I$, and $(f_n')$ converges to $f'$ except possibly at countably many points of $I$.
\end{Theorem}
Observe that the above theorem works if one replaces words "convex" by "concave". Now, by Lemma \ref{concaveEqv} and Theorem \ref{convex_conv}, we conclude this section by the following result.

\begin{Lemma}
\label{uniform}
Let $\Psi$ be an absolutely continuous function on each closed subinterval of $I$ with finite and non zero derivative $\Psi'$ a.e. on $I$.
 If $f_n:I\to\R$ is a sequence of $\Psi$-concave functions converging to a function $f$ which is finite on $I$, then $f$ is $\Psi$-concave on $I$ and the convergence is uniform on any compact subinterval of $I$. Moreover, $(f_n')$ converges to $f'$ a.e. on $I$.
\end{Lemma}

\section{Description of the dual space}


Let $I=(0,l)$, $0<l\leqslant\infty$ and $0<w\in L_0(I)$.
The \emph{weighted Ces{\`a}ro function space} on $I$ is defined to be ($1\leqslant p<\infty$)
\[
C_{p,w}=C_{p,w}(I):=\left\{f\in L_0(I):\|f\|_{C_{p,w}}:=\left(\int_I\left(w(x)\int_0^x|f(t)|\, dt\right)^p\, dx\right)^{1/p}<\infty\right\}\text{.}
\]
Note that for $f\in C_{p,w}(I)$,
\[
\|f\|_{C_{p,w}}=\|\A_w f\|_p\ \ \text{ where }\ \ \A_w f(x)=w(x)\int_0^x|f(t)|\, dt\text{, } \ \ x\in I,
\] and $\|\cdot\|_p$ is the norm in the Lebesgue space $L_p(I)$.

The goal of this section is Theorem \ref{main_th} which gives an isometric description of the Banach dual space $(C_{p,w})^{*}$. We start with two basic lemmas.
\begin{Lemma}
The space $(C_{p,w},\|\cdot\|_{C_{p,w}})$ is an order continuous Banach function lattice with the Fatou property.
\begin{proof}
To see that $(C_{p,w},\|\cdot\|_{C_{p,w}})$ has the Fatou property it is enough to apply Fatou's Lemma twice. Using the Monotone Convergence Theorem one can show that $C_{p,w}(I)$ is an order continuous space \cites{KPS, MR928802}.
\end{proof}
\end{Lemma}


\begin{Lemma}
\label{notinL1}
\begin{enumerate}[{\rm(a)}]
\item The space $C_{p,w}(I)\ne\{0\}$ if and only if $\int_{c}^{l}w(x)^p\, dx<\infty$ for some $c\in I$.
\item $C_{p,w}(I)$ is not continuously embedded into $L_1(I)$ whenever it is not trivial.
\end{enumerate}
\begin{proof}
\begin{enumerate}[(a)]
\item Suppose that $\int_{c}^{l}w(x)^p\, dx<\infty$ for some $c\in I$. For all $d\in(c,l)$ we have that 
\[
\|\chi_{(c,d)}\|_{C_{p,w}}^p=\int_c^l(w(x)\int_0^x\chi_{(c,d)}(t)\, dt)^p\, dx\leqslant (d-c)^p\int_c^lw(x)^p\, dx<\infty\text{,}
\]
 whence $\chi_{(c,d)}\in C_{p,w}$. If $C_{p,w}(I)\ne\{0\}$ then $\chi_{(c,d)}\in C_{p,w}$ for some $c,d\in I$, $d>c$. It follows that $\int_d^lw(x)^p\, dx<\infty$. 
\item Let $a_n$ be a strictly increasing sequence in $I$ such that $\int_{a_n}^{l}w(x)^p\, dx=1/n^p$, $n\geqslant n_0$, for some large enough $n_0\in\N$. If $l=\infty$ let $g_n=\chi_{(a_n,a_n+n)}$, if $l<\infty$ let $g_{n}=\frac{n}{a_{n+1}-a_{n}}\chi
_{(a_{n},a_{n+1})}$. Clearly in both cases $\|g_n\|_1\to\infty$ as $n\to\infty$, $\int_0^xg_n(t)\, dt=0$ for $x<a_n$ and $\int_0^xg_n(t)\, dt\leqslant n$ for $x\geqslant a_n$. Hence $\|g_n\|_{C_{p,w}}^p\leqslant \int_{a_n}^ln^pw(x)^p\, dx=1$ for all $n\geqslant n_0$, and the claim follows.
\end{enumerate}
\end{proof}
\end{Lemma}

If $p=1$ then by Fubini's Theorem
\[
\int_0^{l}w(x)\int_0^x|f(t)|\, dt\, dx=\int_0^l|f(t)|\int_t^lw(x)\, dx\, dt\text{.}
\]
Hence the space $C_{1,w}(I)$ is just a weighted Lebesgue space with weight $\int_t^lw(x)\, dx$, $t\in I$.

In the sequel we assume that $1< p<\infty$ is fixed and the weight function $w$ satisfies the following conditions
\begin{enumerate}[(i)]
\item \label{i}$w>0$ a.e. on $I$,
\item \label{ii}$\int_x^{l}w(t)^p\, dt<\infty$ for all $x\in I$,
\item \label{iii}$\int_0^{l}w(t)^p\, dt=\infty$.
\end{enumerate}
Let further
\begin{equation*}
\Psi(x)=\int_x^{l}w(t)^p\, dt\text{,}\quad x\in I\text{.}
\end{equation*}
Conditions (\ref{i})-(\ref{iii}) imply that the function $\Psi$
is strictly decreasing on $I$, $\lim_{x\to 0^{+}}\Psi(x)=\infty$ and $\lim_{x\to l}\Psi(x)=0$. Also by absolute continuity of $\Psi$ on each compact subinterval of $I$, $\Psi'=-w^p < 0$ a.e. on $I$. Moreover, if $f\in L_0(I)$ is such that $\hat{f}<\infty$ on $I$ then by definition of $D_{\Psi}^{+}\hat{f}$, we get that
\[
D_{\Psi}^{+}\hat{f}(x)=-\hat{f}{'}(x)/w(x)^p\quad\text{for a.a.}\ x\in I\text{,}
\]
where $\hat{f}{'}(x)$ denotes the derivative of $\hat{f}$ at $x$. Note that this derivative exists a.e. because $\Psi$ is absolutely continuous on every closed subinterval of $I$, and so is $\hat{f}$ by Lemma \ref{continuity}.

It is easy to check that if the weight $w$ is a power function $w(x)=x^s$ then conditions (\ref{i})-(\ref{iii}) are satisfied for $s<-1/p$ if $l=\infty$, and for $s\leqslant-1/p$ if $l<\infty$. If $s=-1$ then the space $C_{p,w}$ is the standard Ces{\`a}ro function space $Ces_p$ considered by several authors (see \cite{MR2639977} and the references given therein).

For any $1<p<\infty$ we define its conjugate exponent $q=p/(p-1)$. Let us denote
\begin{equation*}
\B_w f(x)=-\hat{f}{'}(x)/w(x)\quad\text{for a.a. } x\in I\text{.}
\end{equation*}
 We will show that the K\"{o}the dual space of $C_{p,w}(I)$,
\[
(C_{p,w}){'}=(C_{p,w}(I)){'}=\left\{f\in L_0(I):\int_If(t)g(t)\, dt<\infty\text{ for all } g\in C_{p,w}\right\}\text{,}
\]
equipped with the usual norm
\[
\|f\|_{(C_{p,w}){'}}=\sup\left\{\int_If(t)g(t)\, dt:g\in C_{p,w},\|g\|_{C_{p,w}}\leqslant 1\right\}\text{,}
\]
is the space
\[
(C_{p,w}(I)){'}=\left\{f\in L_0(I): \hat{f}<\infty\text{ on } I,\quad\lim_{x\to l}\hat{f}(x)=0\quad\text{and } \B_wf\in L_q(I)\right\}\text{,}
\]
where $\|f\|_{(C_{p,w}){'}}=\|\B_wf\|_q$ and the essential $\Psi$-concave majorant $\hat{f}$ is obtained with respect to $\Psi$. Observe that $\Psi$ (and hence $\hat{f}$) depends on both $p$ and $w$. Since $C_{p,w}$ is an order continuous space with the Fatou property its K\"{o}the dual $(C_{p,w}){'}$ can be identified with its Banach dual space $(C_{p,w})^{*}$. In fact each bounded linear functional $F\in (C_{p,w})^{*}$ is of the integral form $F(g)=\int_If(t)g(t)\, dt$, $g\in C_{p,w}$, where $f\in (C_{p,w}){'}$ and $\|F\|_{(C_{p,w})^{*}}=\|f\|_{(C_{p,w}){'}}$ \cite{KPS, MR928802}.

We start with several preparatory lemmas.
\begin{Lemma}
\label{finitehat2}
If $0\leqslant f\in (C_{p,w}){'}$ then $\esssup_{x\in(y,l)}f(x)<\infty$ and $\esssup_{x\in(0,y)}\frac{f(x)}{\Psi(x)}<\infty$ for all $y\in I$. Consequently $\hat{f}(y)<\infty$ for all $y\in I$.
\begin{proof}

Let $y\in I$ be fixed and $0\leqslant f\in (C_{p,w}){'}$. Suppose that $\esssup_{x\in(y,l)}f(x)=\infty$. For all $C>0$ there exists set $A\subset(y,l)$ with $0<mA<\infty$ such that $f>C$ on $A$. Letting $g=\chi_A/mA$, it follows that $\int_If(t)g(t)\, dt\geqslant C$. But for all sets $A\subset(y,l)$ with positive measure $\|\chi_A/mA\|_{C_{p,w}}\leqslant\Psi(y)^p$, hence $f\notin(C_{p,w}){'}$.

Now we show that $\int_0^y(w(x)/\Psi(x))^p\, dx<\infty$ for all $y\in I$. 
 Fix $y\in I$. By (\ref{ii}) and (\ref{iii}) we can find a sequence $(a_n)$ decreasing to $0$, $a_0=l$ such that 
\[
n\leqslant\int_{a_{n+1}}^{a_n}w(x)^p\, dx<n+1\text{, }n=0,1, \ldots\text{.}
\]
Hence, for $x\in[a_{n+1},a_n)$, $n=0, 1, \ldots$ 
\[
\Psi(x)\geqslant\int_{a_n}^{l}w(t)^p\, dt=\sum_{i=0}^{n-1}\int_{a_{i+1}}^{a_i}w(t)^p\, dt\geqslant\sum_{i=0}^{n-1}i=\frac{n(n-1)}{2}\text{.}
\]
Since $y\in[a_{k+1},a_k)$ for some $k = 0, 1, \ldots$, and $p>1$ we get that
\[
\int_0^y\left(\frac{w(x)}{\Psi(x)}\right)^p\, dx\leqslant
\sum_{n=k}^{\infty}\int_{a_{n+1}}^{a_n}\left(\frac{w(x)}{\Psi(x)}\right)^p\, dx\leqslant\sum_{n=k}^{\infty}\left(\frac{2}{n(n-1)}\right)^p\int_{a_{n+1}}^{a_n}w(x)^p\, dx\leqslant \sum_{n=k}^{\infty}\frac{2^p(n+1)}{n^p(n-1)^p}<\infty\text{.}
\]
Next, for arbitrary $y\in I$, since $\int_0^y(w(x)/\Psi(x))^p\, dx<\infty$, $\Psi$ is decreasing and $1/\Psi$ is bounded on $(0,y)$, denoting $B=\int_0^y1/\Psi(x)\, dx$, by (\ref{ii}) we get that
\begin{equation*}
\begin{split}
\int_I\left(w(x)\int_0^x\frac{\chi_A(t)}{\Psi(t)mA}\, dt\right)^p\, dx&\leqslant\int_0^y\left(\frac{w(x)}{\Psi(x)}\int_0^x\frac{\chi_A(t)}{mA}\, dt\right)^p\, dx+\int_y^l\left(Bw(x)\right)^p\, dx\\
&\leqslant\int_0^y\left(\frac{w(x)}{\Psi(x)}\right)^p\, dx+B^{p}\int_y^lw(x)^p\, dx<\infty\text{.}
\end{split}
\end{equation*}
Hence there is a constant $E>0$ such that $\|1/((mA)\Psi)\chi_A\|_{C_{p,w}}\leqslant E<\infty$ for all $A\subset (0,y)$ with $mA>0$.

Suppose now that $\esssup_{x\in(0,y)}\frac{f(x)}{\Psi(x)}=\infty$. Then for every $C>0$ there exists a set $A\subset (0,y)$ with $mA>0$ such that $f(x)\geqslant C\Psi(x)$ for $x\in A$. Let $g=1/((mA)\Psi)\chi_A$. It follows that $\int_If(x)g(x)\, dx\geqslant\int_AC\Psi(x)(1/((mA)\Psi(x))\, dx=C$ and hence $f\notin(C_{p,w}){'}$ which gives a contradiction.

Finally by Lemma \ref{finitehat1} we have that $\hat{f}(y)<\infty$ for all $y\in I$, and the proof is completed.

\end{proof}
\end{Lemma}

\begin{Lemma}
\label{limf0}
If $f\in L_0(I)$ is such that $\hat{f}\in(C_{p,w}){'}$ then $\lim_{x\to l}\hat{f}(x)=0$.
\begin{proof}
By Corollary \ref{fdecr},  $\hat{f}$ is decreasing and hence $\lim_{x\to l}\hat{f}(x)$ exists.
Suppose that $\lim_{x\to l}\hat{f}(x)=C$ for some $C>0$. It follows that $\hat{f}(t)\geqslant C$ on $I$. By Lemma \ref{notinL1}(b) there is a
sequence of functions $g_n\in B_{C_{p,w}}$ such that $\|g_n\|_1\to\infty$ as $n\to\infty$. Therefore $\int_I\hat{f}(t)g_n(t)\, dt\geqslant C\int_Ig_n(t)\, dt\to\infty\text{,}$ and so $\hat{f}\notin (C_{p,w}){'}$.
\end{proof}

\end{Lemma}

\begin{Lemma}
\label{dpsi1g}
If $f\in L_0(I)$ is such that $\B_wf\in L_q(I)$  then $\lim_{x\to 0^{+}}D_\Psi^{+}\hat{f}(x)=0$.
\begin{proof}
By Proposition \ref{Dpsi} function $D_\Psi^{+}\hat{f}$ is increasing and so $\lim_{x\to 0^{+}}D_\Psi^{+}\hat{f}(x)$ exists. Moreover, (\ref{iii}) and Lemma \ref{Psiinfty} imply that this limit is nonnegative. Suppose that $\lim_{x\to 0^{+}}D_\Psi^{+}\hat{f}(x)=C$ for some constant $C>0$.
Since $D_\Psi^{+}\hat{f}(x)=-\hat{f}{'}(x)/w(x)^p$ a.e., $-\hat{f}{'}(x)/w(x)^p\geqslant C>0$ a.e. It follows $(-\hat{f}{'}(x)/w(x))^q\geqslant C^qw(x)^p$ a.e., so by (\ref{iii}) the integral $\int_I(-\hat{f}{'}(x)/w(x))^q\, dx=\int_I(\B_wf)^q(x)\, dx$ diverges, and this is a contradiction.
\end{proof}
\end{Lemma}

\begin{Lemma}
\label{lem1g}
Let $f\in L_0(I)$ with $\hat{f}<\infty$ on $I$ be such that $\lim_{x\to l}\hat{f}(x)=0$. Then $\int_If(t)g(t)\, dt\leqslant\|\B_wf\|_q\|g\|_{C_{p,w}}$ for any $g\in C_{p,w}$. Consequently $\|f\|_{(C_{p,w}){'}}\leqslant\|\B_wf\|_q$.
\begin{proof}
If $\|\B_wf\|_q=\infty$ then the claim is clear. Assume that $\|\B_wf\|_q<\infty$. In view of $\lim_{x\to l}\hat{f}(x)=0$, using Fubini's Theorem and the H\"{o}lder inequality, for any $g\in C_{p,w}$ we have the following
\begin{equation*}
\begin{split}
\int_0^lf(t)g(t)\, dt&\leqslant\int_0^l|f(t)||g(t)|\, dt\leqslant\int_0^l\hat{f}(t)|g(t)|\, dt=\int_0^l\int_t^l-\hat{f}{'}(x)\, dx |g(t)|\, dt\\
&=\int_0^l\frac{-\hat{f}{'}(x)}{w(x)}w(x)\int_0^x|g(t)|\, dt\, dx\leqslant\|\B_wf\|_q\left\|\A_wg\right\|_p=\|\B_wf\|_q\|g\|_{C_{p,w}}\text{.}
\end{split}
\end{equation*}
From the above it follows that $f\in (C_{p,w}){'}$ and $\|f\|_{(C_{p,w}){'}}\leqslant\|\B_wf\|_q$.
\end{proof}
\end{Lemma}

\begin{Theorem}
\label{thm_psicg}
If $f\in L_0(I)$ is such that $\|\B_wf\|_q<\infty$ and  $\lim_{x\to l}\hat{f}(x)=0$ then $f\in (C_{p,w}){'}$ and $\|f\|_{(C_{p,w}){'}}=\|\B_wf\|_q$.
\begin{proof}

Without loss of generality we assume that $\|\B_wf\|_q=1$. By Lemma \ref{lem1g} we have that $\int_If(t)g(t)\, dt\leqslant\|\B_wf\|_q\|g\|_{C_{p,w}}$, whence $f\in(C_{p,w}){'}$ and $\|f\|_{(C_{p,w}){'}}\leqslant 1$.

The assumption $\|\B_wf\|_q=(\int_I(-\hat{f}{'}(x)/w(x))^q\, dx)^{1/q}<\infty$ gives that $\hat{f}<\infty$ on $I$.
Let $h=(D_\Psi^{+}\hat{f})^{q/p}$. Function $h$ is increasing, finite and right continuous on $I$ and $\lim_{x\to 0^{+}}h(x)=0$ as shown in Proposition \ref{Dpsi} and Lemma \ref{dpsi1g}. Note that $h=(-\hat{f}{'}/w^p)^{q/p}$ a.e. on $I$.

Fix $\epsilon\in(0,1)$. Our main goal is to define a function $g$ such that $\|g\|_{C_{p,w}}\leqslant 1+\epsilon$ and $\int_If(t)g(t)\, dt>1-2\epsilon$. The construction of $g$ will involve a special set $A\subset I$ and a carefully chosen subdivision of $I$. First we find $A$ and then a finite sequence $(a_n)\subset I$ which divides $I$.

Given $y\in I$ let
\[
A_y=\{x\in I:\hat{f}(x)\leqslant |f(x)|+\epsilon/4h(y)\}\text{.}
\]
  Suppose $m(A_y\cap(0,y))=0$ for all $y>y_0$ and some $y_0\in I$. 
 In such case, Lemma \ref{affine_on_interval} implies that $\hat{f}$ is $\Psi$-affine on each interval $(0,y)$, $y\in I$, and hence $\hat{f}$ is $\Psi$-affine on $I$.
By Lemmas \ref{psiaff2} and \ref{dpsi1g}, $D_{\Psi}^{+}\hat{f}$ is constant on $I$ and $\lim_{x\to 0^{-}} D_{\Psi}^{+}\hat{f}(x)=0$. Hence $D_{\Psi}^{+}\hat{f}=0$ on $I$ and so $\hat{f}=0$ by $\lim_{x\to l}\hat{f}(x)=0$, which gives a contradiction with the condition $\|\B_wf\|_q=1$.
This shows that for all $y\in I$ there is $b\in(y,l)$ such that $m(A_b\cap(0,b))>0$.
Since $L_q(I)$ is order continuous we choose $b\in I$ such that
\begin{equation*}
\|(\B_wf)\chi_{(b,l)}\|_q^q\leqslant\epsilon^p/2\text{, and } m(A\cap(0,b))>0\text{,  where }A=A_b\text{,}
\end{equation*}
and $b$ is a point of continuity of $h$.

Observe now, that if there is $x\in I$ such that $m(A\cap(0,x))=0$ then, by Lemma \ref{affine_on_interval}, $\hat{f}$ is $\Psi$-affine on $(0,x)$ and so $h=0$ on $(0,x)$.
Suppose now that $m(A\cap(0,x))>0$ for all $x\in I$. Then for all $y\in I$, there exists $a<y$ such that $m(A\cap(a,x))>0$ whenever $x>a$. Indeed, otherwise there exists $y\in I$ such that for all $a<y$, $m(A\cap(a,x))=0$ for $x>a$ close enough to $a$, say for $x<x_a$. Now the family of sets $(a,x_a)$, $a\in(0,y)$, covers $[\eta,y-\eta]$ for any $\eta>0$. From this, using compactness, one can infer that $m(A\cap(0,y))=0$, which gives a contradiction.

Hence we fix $a\in I$, without loss of generality a point of continuity of $h$, such that $\|(\B_wf)\chi_{(0,a)}\|_q^q\leqslant\epsilon^p/2$ and either 
\begin{equation}
\label{a0nz}
m(A\cap(a,x))>0\text{ for all }x>a\text{, }
\end{equation}
or
\begin{equation}
\label{a0z}
h(a)=0\text{ and }m(A\cap(a,c))=0\text{ for some }c>a\text{.}
\end{equation}
It follows that $b>a$, $m(A\cap(a,b))>0$ and
\begin{equation}
\label{smallpart}
\|(\B_wf)\chi_{(0,a)\cup(b,l)}\|_q^q\leqslant\epsilon^p
\end{equation}

Let $\gamma=\Psi(a)^{1/p}$ and $y_i$ be points of discontinuity of $D_{\Psi}^{+}\hat{f}$ (and hence of $h$) in $(a,b)$ such that $h(y_i^{+})-h(y_i^{-})\geqslant\epsilon/4\gamma$. Here $h(y_i^{+})=\lim_{x\to y_i^{+}}h(x)$ and $h(y_i^{-})=\lim_{x\to y_i^{-}}h(x)$. Clearly there is only finite number of them, say $a<y_1<y_2<\ldots<y_M<b$.
Since $\hat{f}$ is continuous on $I$ for each $y_i$, $i=1, 2,\ldots, M$, we can find two points of continuity of $h$, $\underline{y}_i$,  $\overline{y}_i\in(a,b)$ such that $\underline{y}_i<y_i<\overline{y}_i$, the intervals $[\underline{y}_i,\overline{y}_i]$ are pairwise disjoint,
\begin{equation}
\label{psismallg}
\int_{\underline{y}_i}^{\overline{y_i}}w(x)^p\, dx\leqslant\frac{\epsilon^p}{2^pM(h(y_i^{+})-h(y_i^{-}))^p}\text{,}
\end{equation}


\begin{equation}
\label{fhatclose}
\hat{f}(\underline{y}_i)-\hat{f}(\overline{y}_i)\leqslant\frac{\epsilon}{4h(b)},
\end{equation}

\begin{equation}
\label{hd1g}
h(x)-h(y_i^{+})\leqslant\frac{\epsilon}{4\gamma}\quad\text{for }x\in (y_i,\overline{y}_i),\ \ \text{and}
\end{equation}
\begin{equation}
\label{hd2g}
h(y_i^{-})-h(x)\leqslant\frac{\epsilon}{4\gamma}\quad\text{for }x\in (\underline{y}_i,y_i)\text{.}
\end{equation}

By Lemma \ref{affine_on_interval}, $m(A\cap(\underline{y}_i,\overline{y}_i))>0$ for all $i=1,2,\ldots, M$. Condition (\ref{fhatclose}) implies that
\begin{equation}
\label{hd3g}
\frac{1}{m(A\cap(\underline{y}_i,\overline{y}_i))}\int_{A\cap(\underline{y}_i,\overline{y}_i)}\hat{f}(t)\, dt\geqslant\hat{f}(\underline{y}_i)-\frac{\epsilon}{4h(b)} \text{.}
\end{equation}

Now, the set $(a,b)\setminus\cup_{i=1}^{M}[\underline{y}_i,\overline{y}_i]$ is a union of finite number of open disjoint intervals, say $\cup_{j}(\underline{v}_j,\overline{v}_j)$. Each such interval $(\underline{v}_j,\overline{v}_j)$ can be divided using finite number of points of continuity of $h$, say $u_k$, into subintervals $(u_{k},u_{k+1})$ in such a way that the family $(u_k,u_{k+1})_k$ is a partition of $(\underline{v}_j,\overline{v}_j)$, and
\begin{equation}
\label{hcg}
h(u_{k+1})-h(u_k)\leqslant\frac{\epsilon}{4\gamma}\text{, }
\end{equation}
\begin{equation}
\label{fhatclose2}
\hat{f}(u_k)-\hat{f}(u_{k+1})\leqslant\frac{\epsilon}{4h(b)}\text{.}
\end{equation}

If $m(A\cap(u_k,u_{k+1}))>0$ then by (\ref{fhatclose2}),
\begin{equation}
\label{hc3g}
\frac{1}{m(A\cap(u_k,u_{k+1}))}\int_{A\cap(u_k,u_{k+1})}\hat{f}(t)\, dt\geqslant\hat{f}(u_k)-\frac{\epsilon}{4h(b)} \text{.}
\end{equation}

Let $(a_n)_{n=0}^{N+1}$ be a strictly increasing sequence consisting of all points $u_k\in\cup_{j}(\underline{v}_j,\overline{v}_j)$ obtained above, points $a$, $b$ and $\underline{y}_i$, $\overline{y}_i$, $i=1, 2,\ldots, M$. Note that $a_0=a$, $a_{N+1}=b$ and each interval $(a_n,a_{n+1})$ contains at most one of the points $y_i$, $i=1, 2,\ldots, M$. Note that the union of all sets $(a_n,a_{n+1}]$ is $(a,b]$.

Denote $A_n=A\cap(a_n,a_{n+1})$, $n=0,1,\ldots, N$. Let $E=\{n\in\{0, 1,\ldots,  N\}:mA_n>0\}$. Clearly $E\ne\emptyset$. We can write $E=\{n_1,n_2,\ldots, n_k\}$ where $0\leqslant n_1<n_2<\ldots<n_k\leqslant N$. Note that, if $n\notin E$ then by Lemmas \ref{affine_on_interval} and \ref{psiaff2}, $h(a_{n+1})-h(a_n)=0$ since $h$ is 
constant on $(a_n,a_{n+1})$ and continuous at each point $a_{i}$, $i=0, 1, \ldots N+1$.

Let $\kappa=0$ if $n_1>0$, i.e. if $mA_0=0$ (which is possible when (\ref{a0z}) holds true), $\kappa=h(a_0)/mA_{0}$ if $n_1=0$, i.e. if $mA_0>0$ (which is always a case when (\ref{a0nz}) holds true). Define function
\begin{equation*}
g=\left(\sum_{i=1}^k\frac{h(a_{n_i+1})-h(a_{n_i})}{mA_{n_i}}\chi_{A_{n_i}}+\kappa\chi_{A_0}\right)\sign f\text{.}
\end{equation*}
Now we show that
\[
\|g\|_{C_{p,w}}\leqslant 1+ \epsilon\text{.}
\]
It is clear that $\int_0^x|g(t)|\, dt=0$ if $x<a_0$. Since $h$ is increasing we get that $\int_0^x|g(t)|\, dt\leqslant h(a_{N+1})$ if $x\geqslant a_{N+1}$ and $\int_0^x|g(t)|\, dt\leqslant h(a_{n+1})$ if $a_n\leqslant x<a_{n+1}$, $n=0, 1, \ldots, N$. If $x\in(a_{n},a_{n+1})$ and $(a_{n},a_{n+1})$ does not contain any of points $y_i$, $i=1,2,\dots, M$, then $h(a_{n+1})-h(x)\leqslant\epsilon/4\gamma$ by (\ref{hcg}). Similarly, if $x\in(a_{n},a_{n+1})$ and $(a_{n},a_{n+1})$ contains point $y_i$ for some $i=1,2,\dots, M$, then $h(a_{n+1})-h(x)\leqslant\epsilon/2\gamma+(h(y_i^{+})-h(y_i^{-}))$ by (\ref{hd1g}) and (\ref{hd2g}).
It follows that for $x\in I$,
\[
\A_wg(x)=w(x)\int_0^x|g(t)|\, dt\leqslant w(x)h(x)+\frac{\epsilon}{2\gamma}w(x)\chi_{[a,l)}(x)+w(x)\sum_{i=1}^{M} (h(y_{i}^{+})-h(y_{i}^{-}))\chi_{[\underline{y}_i,\overline{y}_i]}(x)\text{.}
\]

By the triangle inequality, definition of $\gamma$ and  (\ref{psismallg}) we get that
\[
\left\|\frac{\epsilon}{2\gamma}w\chi_{(a,l)}+w\sum_{i=1}^{M} (h(y_{i}^{+})-h(y_{i}^{-}))\chi_{[\underline{y}_i,\overline{y}_i]}\right\|_p\leqslant\epsilon\text{.}
\]
Moreover, $\|wh\|_p=\left(\int_Iw(x)^p(-\hat{f}{'}(x)/w(x)^p)^q\, dx\right)^{1/p}=\left(\|\B_wf\|_q^q\right)^{1/p}=1$ and hence
\[
\|\A_w g\|_p\leqslant\|wh\|_p+\left\|\frac{\epsilon}{2\gamma}w\chi_{(a,l)}+w\sum_{i=1}^{M} (h(y_{i}^{+})-h(y_{i}^{-}))\chi_{[\underline{y}_i,\overline{y}_i]}\right\|_p\leqslant 1+\epsilon\text{.}
\]
Since $\int_A|g(t)|\epsilon/4h(b)\, dt\leqslant\epsilon\int_0^b|g(t)|\, dt/4h(b)\leqslant\epsilon/4$, by definition of $A$ we have that
\begin{equation*}
\begin{split}
\int_If(t)g(t)\, dt&=\int_I|f(t)|g(t)\sign f(t)\, dt\\
&\geqslant\int_A(\hat{f}(t)-\epsilon/4h(b))g(t)\sign f(t)\, dt\geqslant\int_A\hat{f}(t)g(t)\sign f(t)\, dt-\epsilon/4\text{.}
\end{split}
\end{equation*}
Now by (\ref{hd3g}) and (\ref{hc3g}),
\begin{equation*}
\begin{split}
\int_A\hat{f}(t)g(t)\sign f(t)\, dt&=\sum_{i=1}^k\frac{h(a_{n_i+1})-h(a_{n_i})}{mA_{n_i}}\int_{A_{n_i}}\hat{f}(t)\, dt+\kappa\int_{A_0}\hat{f}(t)\, dt\\
&\geqslant \sum_{i=1}^k(h(a_{n_i+1})-h(a_{n_i}))(\hat{f}(a_{n_i})-\epsilon/4h(b))+h(a_0)(\hat{f}(a_0)-\epsilon/4h(b))\text{.}
\end{split}
\end{equation*}
It is easy to see that
\[
\left( \sum_{i=1}^k(h(a_{n_i+1})-(a_{n_i}))+h(a_0) \right)\epsilon/4h(b)\leqslant\epsilon/4\text{.}
\]
Since $h(a_0)=0$ if $n_1>0$, $h(a_{n+1})-h(a_n)=0$ if $n\notin E$ and by (\ref{smallpart}),
\begin{equation*}
\begin{split}
&\sum_{i=1}^k(h(a_{n_i+1})-h(a_{n_i}))\hat{f}(a_{n_i})+h(a_0)\hat{f}(a_0)=\sum_{n=0}^{N}(h(a_{n+1})-h(a_n))\hat{f}(a_n)+h(a_0)\hat{f}(a_0)\\
&=h(a_{N+1})\hat{f}(a_{N+1})+\sum_{n=1}^{N}h(a_{n})(\hat{f}(a_{n-1})-\hat{f}(a_{n}))\\
&=\sum_{n=1}^{N}h(a_{n})\int_{a_{n-1}}^{a_{n}}(-\hat{f}{'}(t))\, dt+h(a_{N+1})\int_{a_{N+1}}^l(-\hat{f}{'}(t))\, dt\\
&\geqslant\sum_{n=1}^{N}\int_{a_{n-1}}^{a_{n}}h(t)(-\hat{f}{'}(t))\, dt+\int_{a_{N+1}}^lh(t)(-\hat{f}{'}(t))\, dt \geqslant\int_a^b\left(\frac{-\hat{f}{'}(t)}{w(t)}\right)^q\, dt\geqslant1-\epsilon^p\geqslant 1-\epsilon\text{.}
\end{split}
\end{equation*}
Combining all the above together we obtain $\int_If(t)g(t)\, dt\geqslant 1-3\epsilon/2$. Dividing both sides by $1+\epsilon$ one gets
\[
\int_If(t)g(t)/(1+\epsilon)\, dt\geqslant 1-3\epsilon\text{.}
\]
Finally, by $\|g/(1+\epsilon)\|_{C_{p,w}}\leqslant 1$ it follows that $\|f\|_{(C_{p,w}){'}}=1$.

\end{proof}
\end{Theorem}

\begin{Lemma}
\label{fhat}
If $f\in (C_{p,w}){'}$ then $\hat{f}\in (C_{p,w}){'}$ and $\|f\|_{(C_{p,w}){'}}=\|\hat{f}\|_{(C_{p,w}){'}}=\|\B_wf\|_q$.
\begin{proof}
Consider first the case when $l<\infty$. Let $f\in (C_{p,w}){'}$ and  $f_m=f\chi_{[1/m,l-1/m]}$, $m\in\N$. By Lemma \ref{finitehat2}, $\hat{f}<\infty$ on $I$. Clearly $\widehat{f_m}\leqslant \hat{f}$. Letting $y\in I$, by definition of $\hat{f}$, for every $\epsilon>0$ there exist $n\in\N$ and a set $A=\{(y_1,  \ldots,y_n)\in I^n:\sum_{i=1}^n\alpha_if(y_i)>\hat{f}(y)-\epsilon, \sum_{i=1}^n\alpha_i\Psi(y_i)=\Psi(y), \sum_{i=1}^n\alpha_i=1, \alpha_i\geqslant 0, i=1,2,\ldots, n\}$ with $m^{(n)}A>0$. Let $r>0$ be such that $1/r<y<l-1/r$ and $m^{(n)}(A\cap(1/r,l-1/r)^n)>0$. Since for all $m>r$, $f=f_m$ on $(1/r,l-1/r)$ it follows that $\widehat{f_m}(y)>\hat{f}(y)-\epsilon$. By arbitrariness of $\epsilon$ we get that $\widehat{f_m}(y)\to\hat{f}(y)$ as $m\to\infty$. By Lemma \ref{uniform}, $\widehat{f_m}^{'}\to\hat{f}{'}$ a.e. on $I$.

Note that $D_{\Psi}^{+}\widehat{f_m}(x)=-\widehat{f_m}'(x)/w(x)^p$ is $0$ a.e. on $(0,1/m)$ and constant a.e. on $(l-1/m,l)$. Hence the function $-\widehat{f_m}'(x)/w(x)=w(x)^{p-1}D_{\Psi}^{+}\widehat{f_m}(x)$ a.e. is
in $L_q(I)$, $m\in\N$. By Theorem \ref{thm_psicg} we get that $\|f_m\|_{(C_{p,w}){'}}=\|\B_wf_m\|_q$ for all $m\in\N$.
Now by Lemma \ref{lem1g}, Lemma \ref{uniform}, the Fatou Lemma and the Fatou property of $(C_{p,w}){'}$ we get that
\begin{equation*}
\begin{split}
\|\hat{f}\|_{(C_{p,w}){'}}&\leqslant\|\B_wf\|_q=\left(\int_I(-\hat{f}'(x)/w(x))^q\, dx\right)^{1/q}=\left(\int_I(-\lim_n\widehat{f_n}'(x)/w(x))^q\, dx\right)^{1/q}\\
&\leqslant\liminf_n\left(\int_I(-\widehat{f_n}'(x)/w(x))^q\, dx\right)^{1/q}
\leqslant\sup_n\|\B_wf_n\|_q=\sup_n\|f_n\|_{(C_{p,w}){'}}=\|f\|_{(C_{p,w}){'}}<\infty\text{.}
\end{split}
\end{equation*}
So $\|f\|_{(C_{p,w}){'}}=\|\hat{f}\|_{(C_{p,w}){'}}$, since $|f|\leqslant\hat{f}$. The above inequality also shows that $\|f\|_{(C_{p,w}){'}}=\|\B_wf\|_q$.

In case when $l=\infty$ we proceed similarly as above taking $f_m=f\chi_{[1/m,m]}$, $m\in\N$.
\end{proof}
\end{Lemma}

Now we are ready to present the main result in this section, isometric description of the dual space $(C_{p,w})^{*}$. Namely, by Lemma \ref{limf0}, Theorem \ref{thm_psicg} and Lemma \ref{fhat} we get the following theorem.
\begin{Theorem}
\label{main_th}
Let $1<p<\infty$, $q=\frac{p}{p-1}$, $\Psi(x)=\int_x^lw(t)^p\, dt$, $x\in I=(0,l)$, $0<l\leqslant\infty$. Then a function $f\in(C_{p,w}){'}$ if and only if $\hat{f}<\infty$ on $I$, $\lim_{x\to l}\hat{f}(x)=0$ and $\|\B_wf\|_q<\infty$. Moreover
\[
 \|f\|_{(C_{p,w}){'}}=\|\B_wf\|_q\text{ for all }f\in(C_{p,w}){'}\text{.}
\]
The Banach dual space $(C_{p,w})^{*}$ of $C_{p,w}$ is isometrically isomorphic to $(C_{p,w}){'}$ in the sense that every $F\in(C_{p,w})^{*}$ is of the form
\[
F(g)=\int_If(t)g(t)\, dt\text{,}\quad g\in C_{p,w}\text{,}
\]
for a unique $f\in(C_{p,w}){'}$ and $\|F\|_{(C_{p,w})^{*}}=\|f\|_{(C_{p,w}){'}}$.
\end{Theorem}

\section{Diameter of slices of the unit ball}

Let $(X,\|\cdot\|)$ be a Banach space. 
Recall that the set $s(x^*;\eta)=\{x\in B_X:x^*x>1-\eta\}$, where $x^*\in S_{X^*}$, $0<\eta<1$, is called a \emph{slice} of $B_X$. Applying the techniques described in the previous section we show the following result.

\begin{Theorem}
\label{diam2}
Every slice of $B_{C_{p,w}}$ has diameter $2$.
\begin{proof}
Let $f\in S_{(C_{p,w}){'}}$ and $0<\eta<1$ be fixed. We will show that the slice $s(F;\eta)$, where $F\in (C_{p,w})^{*}$ is defined by $F(g)=\int_If(t)g(t)\, dt$, $g\in C_{p,w}$, has diameter $2$.

Let $0<\epsilon<\eta/10$ be arbitrary. Let the function $h$, the set $A$, points $\underline{y}_i$,$y_i$, $\overline{y}_i$, $i = 1, 2, \ldots M$, the sequence $(a_n)_{n=0}^{N+1}$, the set $E$ and $\gamma$ be defined exactly in the same way as in the first part of the proof of Theorem \ref{thm_psicg}.
We have that $\|(\B_wf)\chi_{(0,a)\cup(b,l)}\|_q^q\leqslant\epsilon^p$ and
\begin{equation}
\label{epspart}
\left\|\frac{\epsilon}{\gamma}w\chi_{(a,l)}+2w\sum_{i=1}^{M} (h(y_{i}^{+})-h(y_{i}^{-}))\chi_{[\underline{y}_i,\overline{y}_i]}\right\|_p\leqslant 2\epsilon\text{.}
\end{equation}


Let sets $B$ and $C$ be any measurable sets such that $A=B\cup C$, $B\cap C=\emptyset$, and for all $n\in E$, $m(B\cap(a_n,a_{n+1}))>0$ and $m(C\cap(a_n,a_{n+1}))>0$. Denote $B_n=B\cap(a_n,a_{n+1})$ and $C_n=C\cap(a_n,a_{n+1})$, $n=0,1,\ldots, N$. Let $\kappa_1=\kappa_2=0$ if $n_1>0$ (i.e. if $mA_0=0$), $\kappa_1=h(a_0)/mB_0$, $\kappa_2=h(a_0)/mC_0$ if $n_1=0$ (i.e. if $mA_0>0$). Define functions
\begin{equation*}
g_1=\left(\sum_{i=1}^k\frac{h(a_{n_i+1})-h(a_{n_i})}{mB_{n_i}}\chi_{B_{n_i}}+\kappa_1\chi_{B_0}\right)\sign f\text{,}
\end{equation*}
\begin{equation*}
g_2=\left(\sum_{i=1}^k\frac{h(a_{n_i+1})-h(a_{n_i})}{mC_{n_i}}\chi_{C_{n_i}}+\kappa_2\chi_{C_0}\right)\sign f\text{.}
\end{equation*}


Similarly as in the proof of Theorem \ref{thm_psicg} one can show that $\|g_j\|_{C_{p,w}}\leqslant 1+\epsilon$, $j=1,2$. Since formulas (\ref{fhatclose}) and (\ref{fhatclose2}) do not depend on the set $A$, (\ref{hd3g}) and (\ref{hc3g}) remain true if we replace set $A$ by $B$ or $C$. Hence we also get the estimates $\int_If(t)g_j(t)/(1+\epsilon)\,dt\geqslant 1-3\epsilon>1-\eta$, $j=1,2$. It follows that $g_j/(1+\epsilon)\in s(F;\eta)$, $j=1,2$.


Observe that $\int_0^x|g_j(t)|\, dt=0$ if $x<a_0$. If $x\in(a_n,a_{n+1})$ and $y_i\notin (a_n,a_{n+1})$ for all $i=1, 2, \ldots, M$, then $\int_0^x|g_j(t)|\, dt\geqslant h(a_n)\geqslant h(a_{n+1})-\epsilon/2\gamma$ by (\ref{hcg}). Similarly, if $x\in(a_n,a_{n+1})$ and $y_i\in (a_n,a_{n+1})$ for some $i=1, 2, \ldots, M$, then $\int_0^x|g_j(t)|\, dt\geqslant h(a_n)\geqslant h(a_{n+1})-\epsilon/2\gamma-(h(y_i^{+})-h(y_i^{-}))$ by (\ref{hd1g}) and (\ref{hd2g}). It follows that for all $x\in(a,b)$, $j=1, 2$,
\[
\int_0^x|g_j(t)|\, dt\geqslant h(x)\chi_{(a,b)}(x)-\frac{\epsilon}{2\gamma}\chi_{(a,b)}(x)-\sum_{i=1}^M(h(y_i^{+})-h(y_i^{-}))\chi_{(\underline{y}_i,\overline{y}_i)}(x)\text{.}
\]
Since $g_1$ and $g_2$ have disjoint supports, we get that 
for all $x\in(a,b)$,
\begin{equation*}
\begin{split}
w(x)\int_0^x|g_1(t)-g_2(t)|\, dt&=w(x)\int_0^x|g_1(t)|\, dt+w(x)\int_0^x|g_2(t)|\, dt\\
&\geqslant 2w(x)h(x)\chi_{(a,b)}(x)-\frac{\epsilon}{\gamma}w(x)\chi_{(a,b)}(x)-2w(x)\sum_{i=1}^M(h(y_i^{+})-h(y_i^{-}))\chi_{(\underline{y}_i,\overline{y}_i)}(x)\text{.}
\end{split}
\end{equation*}
It follows that
\begin{equation*}
\begin{split}
\|g_1-g_2\|_{C_{p,w}}&\geqslant\left\|2wh\chi_{(a,b)}-
\left(\frac{\epsilon}{\gamma}w\chi_{(a,b)}+2w\sum_{i=1}^M(h(y_i^{+})-h(y_I^{-}))\chi_{(\underline{y}_i,\overline{y}_i)}\right)\right\|_p\\
&\geqslant\left|2\|wh\chi_{(a,b)}\|_p-\left\|\frac{\epsilon}{\gamma}w\chi_{(a,b)}+2w\sum_{i=1}^M(h(y_i^{+})-h(y_i^{-}))\chi_{(\underline{y}_i,\overline{y}_i)}\right\|_p\right|\text{.}
\end{split}
\end{equation*}
Since $\|wh \|_p=1$ and $\|wh\chi_{(0,a)\cup(b,l)}\|_p^p=\|(\B_wf)\chi_{(0,a)\cup(b,l)}\|_q^q\leqslant \epsilon^p$ we get that
\begin{equation*}
\begin{split}
\|wh\chi_{(a,b)}\|_p&=\|wh-wh\chi_{(0,a)\cup(b,l)} \|_p\geqslant \left|\|wh\|_p-\|wh\chi_{(0,a)\cup(b,l)}\|_p \right|\\
&=1- \|(\B_wf)\chi_{(0,a)\cup(b,l)}\|_q^{q/p}\geqslant 1-\epsilon\text{.}
\end{split}
\end{equation*}
By the above and (\ref{epspart})
we get that $\|g_1-g_2\|_{C_{p,w}}\geqslant 2-4\epsilon$. Dividing now both sides by $1+\epsilon$ we obtain that $\|(g_1-g_2)/(1+\epsilon)\|_{C_{p,w}}\geqslant 2-6\epsilon$. Since $\epsilon$ can be taken arbitrarily small we obtain that the diameter of $s(F;\eta)$ is $2$.

\end{proof}
\end{Theorem}

\section{Final conclusions}
Recall that a Banach space $(X,\|\cdot\|)$ is called \emph{locally uniformly convex} if for any $x\in S_X$ and any sequence $(x_n)\subset B_X$, $\lim_{n\to\infty}\|x+x_n\|=2$ implies that $\lim_{n\to\infty}\|x-x_n\|=0$. A point $x\in S_X$ is said to be \emph{strongly exposed} if there is $x^{*}\in S_{X^{*}}$ such that $x^{*}x=1$, $x^{*}y<1$ for all $y\in B_X\setminus\{x\}$, and $x^{*}x_n\to 1$ implies that $\|x-x_n\|\to 0$ as $n\to\infty$ for any sequence $(x_n)\subset B_X$. 

A point $x\in S_X$ is called a \emph{denting point} of $B_X$ if $x\notin\overline{co}\{B_X\setminus(x+\epsilon B_X)\}$ for each $\epsilon>0$. It is easy to see that if the unit ball $B_X$ has denting points then it has slices of arbitrary small diameter \cite[Proposition 2.3.2, p. 28]{MR704815}. Also any strongly exposed point is a denting point \cite[p. 227]{MR2300779} and in locally uniformly convex space all points of its unit sphere are denting \cite{MR0098976}. The Radon-Nikodym property can be characterized in terms of denting points. Namely, a Banach space $X$ has the Radon-Nikodym property if and only if for every equivalent norm in $X$ the respective unit ball $B_X$ has a denting point \cite[p. 30]{MR704815}. For definition and more details on Radon-Nikodym property we refer to \cite{MR704815}.
Consequently by Theorem \ref{diam2} we get the following corollaries.

\begin{Corollary}
The space $(C_{p,w},\|\cdot\|_{C_{p,w}})$ does not have the Radon-Nikodym property.
\end{Corollary}

\begin{Corollary}
The unit sphere of $(C_{p,w},\|\cdot\|_{C_{p,w}})$ does not have strongly exposed points.
\end{Corollary}

\begin{Corollary}
The unit sphere of $(C_{p,w},\|\cdot\|_{C_{p,w}})$ does not have denting points.
\end{Corollary}

\begin{Corollary}
The space $(C_{p,w},\|\cdot\|_{C_{p,w}})$ is not locally uniformly convex.
\end{Corollary}

\begin{Corollary}
The space $(C_{p,w},\|\cdot\|_{C_{p,w}})$ is not a dual space.
\begin{proof}
It is known that every separable dual space has the Kre\u\i n-Milman Property \cite{MR0211244} and that the latter is equivalent to the Radon-Nikodym Property in Banach lattices \cite{MR589144, MR934876}. Since $C_{p,w}$ is a separable Banach lattice without the Radon-Nikodym Property it cannot be a dual space.
\end{proof}
\end{Corollary}

\section*{Acknowledgment}
The authors thank referees for valuable comments and suggestions.

\bibliography{Ces3}

\def\polhk#1{\setbox0=\hbox{#1}{\ooalign{\hidewidth
  \lower1.5ex\hbox{`}\hidewidth\crcr\unhbox0}}} \def\cprime{$'$}
\begin{bibdiv}
\begin{biblist}

\bib{pr68}{article}{
       title={Programma van {J}aarlijkse {P}rijsvragen ({A}nnual {P}roblem
  {S}ection)},
        date={1968},
     journal={Nieuw Arch. Wisk.},
      volume={16},
       pages={47\ndash 51},
}

\bib{MR2431042}{article}{
      author={Astashkin, Sergei~V.},
      author={Maligranda, Lech},
       title={Ces\`aro function spaces fail the fixed point property},
        date={2008},
        ISSN={0002-9939},
     journal={Proc. Amer. Math. Soc.},
      volume={136},
      number={12},
       pages={4289\ndash 4294},
         url={http://dx.doi.org/10.1090/S0002-9939-08-09599-3},
      review={\MR{2431042 (2009g:46045)}},
}

\bib{MR2639977}{article}{
      author={Astashkin, Sergei~V.},
      author={Maligranda, Lech},
       title={Structure of {C}es\`aro function spaces},
        date={2009},
        ISSN={0019-3577},
     journal={Indag. Math. (N.S.)},
      volume={20},
      number={3},
       pages={329\ndash 379},
         url={http://dx.doi.org/10.1016/S0019-3577(10)00002-9},
      review={\MR{2639977}},
}

\bib{MR2650988}{article}{
      author={Astashkin, Sergei~V.},
      author={Maligranda, Lech},
       title={Rademacher functions in {C}es\`aro type spaces},
        date={2010},
        ISSN={0039-3223},
     journal={Studia Math.},
      volume={198},
      number={3},
       pages={235\ndash 247},
         url={http://dx.doi.org/10.4064/sm198-3-3},
      review={\MR{2650988}},
}

\bib{MR928802}{book}{
      author={Bennett, Colin},
      author={Sharpley, Robert},
       title={Interpolation of operators},
      series={Pure and Applied Mathematics},
   publisher={Academic Press Inc.},
     address={Boston, MA},
        date={1988},
      volume={129},
        ISBN={0-12-088730-4},
      review={\MR{928802 (89e:46001)}},
}

\bib{MR1317938}{article}{
      author={Bennett, Grahame},
       title={Factorizing the classical inequalities},
        date={1996},
        ISSN={0065-9266},
     journal={Mem. Amer. Math. Soc.},
      volume={120},
      number={576},
       pages={viii+130},
      review={\MR{1317938 (96h:26020)}},
}

\bib{MR0211244}{article}{
      author={Bessaga, C.},
      author={Pe{\l}czy{\'n}ski, A.},
       title={On extreme points in separable conjugate spaces},
        date={1966},
        ISSN={0021-2172},
     journal={Israel J. Math.},
      volume={4},
       pages={262\ndash 264},
      review={\MR{0211244 (35 \#2126)}},
}

\bib{MR589144}{article}{
      author={Bourgain, Jean},
      author={Talagrand, Michel},
       title={Dans un espace de {B}anach reticul\'e solide, la propri\'et\'e de
  {R}adon-{N}ikod\'ym et celle de {K}re\u\i n-{M}il\cprime man sont
  \'equivalentes},
        date={1981},
        ISSN={0002-9939},
     journal={Proc. Amer. Math. Soc.},
      volume={81},
      number={1},
       pages={93\ndash 96},
         url={http://dx.doi.org/10.2307/2043994},
      review={\MR{589144 (82d:46032)}},
}

\bib{MR704815}{book}{
      author={Bourgin, Richard~D.},
       title={Geometric aspects of convex sets with the {R}adon-{N}ikod\'ym
  property},
      series={Lecture Notes in Mathematics},
   publisher={Springer-Verlag},
     address={Berlin},
        date={1983},
      volume={993},
        ISBN={3-540-12296-6},
      review={\MR{704815 (85d:46023)}},
}

\bib{MR934876}{article}{
      author={Caselles, V.},
       title={A short proof of the equivalence of {KMP} and {RNP} in {B}anach
  lattices and preduals of von {N}eumann algebras},
        date={1988},
        ISSN={0002-9939},
     journal={Proc. Amer. Math. Soc.},
      volume={102},
      number={4},
       pages={973\ndash 974},
         url={http://dx.doi.org/10.2307/2047343},
      review={\MR{934876 (89c:46029)}},
}

\bib{MR1904283}{incollection}{
      author={Chen, S.},
      author={Cui, Y.},
      author={Hudzik, H.},
      author={Sims, B.},
       title={Geometric properties related to fixed point theory in some
  {B}anach function lattices},
        date={2001},
   booktitle={Handbook of metric fixed point theory},
   publisher={Kluwer Acad. Publ.},
     address={Dordrecht},
       pages={339\ndash 389},
      review={\MR{1904283 (2003f:46031)}},
}

\bib{MR1744077}{article}{
      author={Cui, Yunan},
      author={Hudzik, Henryk},
       title={Some geometric properties related to fixed point theory in
  {C}es\`aro spaces},
        date={1999},
        ISSN={0010-0757},
     journal={Collect. Math.},
      volume={50},
      number={3},
       pages={277\ndash 288},
      review={\MR{1744077 (2001f:46033)}},
}

\bib{MR1972393}{inproceedings}{
      author={Cui, Yunan},
      author={Hudzik, Henryk},
       title={Packing constant for {C}esaro sequence spaces},
        date={2001},
   booktitle={Proceedings of the {T}hird {W}orld {C}ongress of {N}onlinear
  {A}nalysts, {P}art 4 ({C}atania, 2000)},
      volume={47},
       pages={2695\ndash 2702},
         url={http://dx.doi.org/10.1016/S0362-546X(01)00389-3},
      review={\MR{1972393 (2004c:46033)}},
}

\bib{MR1772119}{incollection}{
      author={Cui, Yunan},
      author={Hudzik, Henryk},
      author={Li, Yanhong},
       title={On the {G}arc\'\i a-{F}alset coefficient in some {B}anach
  sequence spaces},
        date={2000},
   booktitle={Function spaces ({P}ozna\'n, 1998)},
      series={Lecture Notes in Pure and Appl. Math.},
      volume={213},
   publisher={Dekker},
     address={New York},
       pages={141\ndash 148},
      review={\MR{1772119 (2001h:46009)}},
}

\bib{MR1810056}{article}{
      author={Cui, Yunan},
      author={Meng, Chenghui},
      author={P{\l}uciennik, Ryszard},
       title={Banach-{S}aks property and property {$(\beta)$} in {C}es\`aro
  sequence spaces},
        date={2000},
        ISSN={0129-2021},
     journal={Southeast Asian Bull. Math.},
      volume={24},
      number={2},
       pages={201\ndash 210},
         url={http://dx.doi.org/10.1007/s100120070003},
      review={\MR{1810056 (2001m:46031)}},
}

\bib{MR1608225}{article}{
      author={Cui, Yunan},
      author={P{\l}uciennik, Ryszard},
       title={Local uniform nonsquareness in {C}es\`aro sequence spaces},
        date={1997},
        ISSN={0373-8299},
     journal={Comment. Math. Prace Mat.},
      volume={37},
       pages={47\ndash 58},
      review={\MR{1608225 (99b:46025)}},
}

\bib{MR0098976}{article}{
      author={Fan, Ky},
      author={Glicksberg, Irving},
       title={Some geometric properties of the spheres in a normed linear
  space},
        date={1958},
        ISSN={0012-7094},
     journal={Duke Math. J.},
      volume={25},
       pages={553\ndash 568},
      review={\MR{0098976 (20 \#5421)}},
}

\bib{MR0348444}{article}{
      author={Jagers, A.~A.},
       title={A note on {C}es\`aro sequence spaces},
        date={1974},
        ISSN={0028-9825},
     journal={Nieuw Arch. Wisk. (3)},
      volume={22},
       pages={113\ndash 124},
      review={\MR{0348444 (50 \#942)}},
}

\bib{KPS}{book}{
      author={Kre{\u\i}n, S.~G.},
      author={Petun{\={\i}}n, Yu.~{\=I}.},
      author={Sem{\"e}nov, E.~M.},
       title={Interpolation of linear operators},
      series={Translations of Mathematical Monographs},
   publisher={American Mathematical Society},
     address={Providence, R.I.},
        date={1982},
      volume={54},
        ISBN={0-8218-4505-7},
        note={Translated from the Russian by J. Sz{\H{u}}cs},
      review={\MR{649411 (84j:46103)}},
}

\bib{MR2178902}{book}{
      author={Niculescu, Constantin~P.},
      author={Persson, Lars-Erik},
       title={Convex functions and their applications},
      series={CMS Books in Mathematics/Ouvrages de Math\'ematiques de la SMC,
  23},
   publisher={Springer},
     address={New York},
        date={2006},
        ISBN={978-0387-24300-9; 0-387-24300-3},
        note={A contemporary approach},
      review={\MR{2178902 (2006m:26001)}},
}

\bib{MR0272960}{article}{
      author={Peetre, J.},
       title={Concave majorants of positive functions},
        date={1970},
        ISSN={0001-5954},
     journal={Acta Math. Acad. Sci. Hungar.},
      volume={21},
       pages={327\ndash 333},
      review={\MR{0272960 (42 \#7841)}},
}

\bib{MR2300779}{book}{
      author={Pietsch, Albrecht},
       title={History of {B}anach spaces and linear operators},
   publisher={Birkh\"auser Boston Inc.},
     address={Boston, MA},
        date={2007},
        ISBN={978-0-8176-4367-6; 0-8176-4367-2},
      review={\MR{2300779 (2008i:46002)}},
}

\bib{MR0442824}{book}{
      author={Roberts, A.~Wayne},
      author={Varberg, Dale~E.},
       title={Convex functions},
   publisher={Academic Press [A subsidiary of Harcourt Brace Jovanovich,
  Publishers], New York-London},
        date={1973},
        note={Pure and Applied Mathematics, Vol. 57},
      review={\MR{0442824 (56 \#1201)}},
}

\bib{MR924157}{book}{
      author={Rudin, Walter},
       title={Real and complex analysis},
     edition={Third},
   publisher={McGraw-Hill Book Co.},
     address={New York},
        date={1987},
        ISBN={0-07-054234-1},
      review={\MR{924157 (88k:00002)}},
}

\bib{MR0276751}{article}{
      author={Shiue, Jau-shyong},
       title={A note on {C}es\`aro function space},
        date={1970},
        ISSN={0049-2930},
     journal={Tamkang J. Math.},
      volume={1},
      number={2},
       pages={91\ndash 95},
      review={\MR{0276751 (43 \#2491)}},
}

\bib{MR2129625}{book}{
      author={Stein, Elias~M.},
      author={Shakarchi, Rami},
       title={Real analysis},
      series={Princeton Lectures in Analysis, III},
   publisher={Princeton University Press},
     address={Princeton, NJ},
        date={2005},
        ISBN={0-691-11386-6},
        note={Measure theory, integration, and Hilbert spaces},
      review={\MR{2129625 (2005k:28024)}},
}

\end{biblist}
\end{bibdiv}



\end{document}